\documentclass[11pt]{amsart}
\usepackage{amssymb}
\usepackage{amsmath}
\usepackage{mathtools, float}
\usepackage{amsfonts}
\usepackage{stmaryrd}
\usepackage[english]{babel}
\usepackage[utf8]{inputenc}
\usepackage{enumerate}
\usepackage{graphicx,pinlabel,xcolor, fullpage}
\usepackage{adjustbox}
\usepackage{array,multirow}
\usepackage{float}

\theoremstyle{plain}
\newtheorem{theorem}{Theorem}[section]
\newtheorem{corollary}[theorem]{Corollary}
\newtheorem{lemma}[theorem]{Lemma}

\newtheorem{theorem*}{Theorem}

\theoremstyle{definition}
\newtheorem{definition}[theorem]{Definition}

\newtheorem{remark}[theorem]{Remark}
\newtheorem{note}[theorem]{Note}
\newtheorem{proposition}[theorem]{Proposition}
\newtheoremstyle{case}{}{}{}{}{}{:}{ }{}
\theoremstyle{case}

\newcommand{\QED}{\hfill $\blacksquare$}
\newcommand{\N}{\mathbb N}

\begin{document}
\title[Finite Symmetries of surfaces of $p$-groups of co-class 1]{Finite Symmetries of surfaces of $p$-groups of co-class 1}
\author{Siddhartha Sarkar}
\address{Department of Mathematics\\
Indian Institute of Science Education and Research Bhopal\\
Bhopal Bypass Road, Bhauri \\
Bhopal 462 066, Madhya Pradesh\\
India}
\email{sidhu@iiserb.ac.in}
\urladdr{https://home.iiserb.ac.in/~sidhu/}



\begin{abstract}
The genus spectrum of a finite group $G$ is a set of integers $g \geq 2$ such that $G$ acts on a closed orientable compact surface $\Sigma_g$ of genus $g$ preserving the orientation. In this paper we complete the study of spectrum sets of finite $p$-groups of co-class $1$, where $p$ is an odd prime. As a consequence we prove that given an order $p^n$ and exponent $p^e$, there are at the most eight genus spectrum despite the infinite growth of their isomorphism types along $(n,e)$. Based on these results we also classify these groups which has unique stable upper genus $\sigma_e(p^e) - p^e$, where $\sigma_e(p)$ is a constant that depends on $p$ and $e$. 
\end{abstract}

\maketitle

\section{Introduction}
\label{introsec}

Let $G$ denote a finite group and $\Sigma_g$ a compact closed orientable surface of genus $g \geq 0$. We say $G$ acts on $\Sigma_g$ whenever there is a Riemann surface structure $(X, \Sigma_g)$ on $\Sigma_g$, so that the action is induced by a faithful and orientation preserving action of $G$ on $X$. The set sp($G$) denote the set of all $g \geq 2$, so that $G$ acts on $\Sigma_g$, and is called the genus spectrum of $G$. 
 
Determining the spectrum set sp($G$) of $G$ is termed as the Hurwitz problem in {\cite{mmi}}. Solving Hurwitz problem for arbitrary finite groups require the knowledge of computing the least stable solution of Frobenius problem with arbitrary large number of variables (see {\cite{owe}}), which is difficult. Hence the spectrum sets are known only for a few finite groups.  

(i) arbitrary finite groups: almost all finite cyclic groups of order $pq$ {\cite{owe}}, metacylic groups of order $pq$, ($p, q$ distinct primes) {\cite{wea}}, some finite simple groups including half of the sporadic simple groups {\cite{ms2}}.

(ii) finite groups of prime power order: cyclic {\cite{kma}}, groups with large cyclic subgroups {\cite{mta}}, exponent $p$ and groups of maximal class of order $\leq p^p$ {\cite{sar}}, certain abelian groups {\cite{ms1}}.

Let exponent of $G$ (denoted by ${\mathrm{exp}}(G)$) denote the least common multiple of the order of the elements of $G$. To understand the spectrum sets of finite $p$-groups, it is reasonable to ask if there is a property $P$ so that for fixed $n, e \in {\mathbb N}$ with $e \leq n$, the groups in the subclass 
\[
{\mathcal {G}}_{n,e,P} = \{ G ~:~ |G|=p^n, {\mathrm{exp}}(G)=p^e, G {\mathrm{~has~property~}} P \}
\]
has only a finite number of spectrum sets. Although there is no complete classification of finite $p$-groups, the most effective way to do so is the co-class theory which is confirmed by the proof of co-class conjecture due to Shalev (see {\cite{sha}}, {\cite{lmc}}).      
 
A good starting point thus should be examining finite $p$-groups of co-class $1$ (i.e., nilpotency class $cl(G)$ of $G$ satisfy $\log_p(|G|) - cl(G) = 1$ and also called $p$-groups of maximal class). For $p=2$, there are only three isomorphism types of such groups in a given fixed order, which has distinct spectrum sets as studied in the PhD thesis of Talu in 1993 (see {\cite{mta}}). For odd primes $p$, it is also known that among the groups of fixed order $\leq p^p$ and exponent there is only one spectrum set {\cite{sar}}. In this paper we prove that if ${\mathcal G}_{n,e,1}$ is the collection of all finite $p$-groups of order $p^n$, exponent $p^e$ and co-class $1$, then the number of spectrum sets of groups $G \in {\mathcal G}_{n,e,1}$ is bounded above by $8$. To see this we thus concentrate to the case when $p$ is and odd prime, $G$ is a $p$-group of co-class $1$ and $|G| \geq p^{p+1}$. 

\smallskip

For a finite $p$-group of maximal class of order $p^n \geq p^{p+1}$ and exponent $p^e$, let $G_2 = [G,G], G_4 = [G,G,G,G]$ and $G_1$ denote the maximal subgroup $G_1 = C_G(G_2/G_4)$ of $G$. For a subset $X \subseteq G$ let $||X||$ denote the set of all order of elements of $X$. Let $G^{\ast}_0 := G \setminus G_1$. For a fixed $|G| \geq p^{p+1}$ and exponent $p^e$, these groups fall into five types: Type-I : $||G^{\ast}_0|| = \{ p \}$, Type-II $||G^{\ast}_0|| = \{ p^2 \}$, and \newline while $||G^{\ast}_0|| = \{ p, p^2 \}$, Type-III : $G$ is not exceptional (require $p \geq 5$), Type-IV : $G$ is $p$-exceptional and Type-V : $G$ is $p^2$-exceptional (see Definition \ref{def-excep}). This make ${\mathcal G}_{n,e,1}$ divided into five subclasses so that every group belong to the same subclass has the same genus spectrum. Using counting methods one can show that the number of groups in each such subclass $\rightarrow \infty$ as $n \rightarrow \infty$. 

\smallskip

To describe the spectrum sets precisely we require another another property of exponent depth (denoted by ${\mathrm{ed}}(G)$): ${\mathrm{ed}}(G) \geq 2$ means that ${\mathrm{exp}}(G_1) = {\mathrm{exp}}(G_2)$, whereas ${\mathrm{ed}}(G) = 1$ means ${\mathrm{exp}}(G_1) > {\mathrm{exp}}(G_2)$. Here we describe the main results of the paper:  

\smallskip

\noindent {\bf Theorem A.} Let $G$ be a finite $p$-group of maximal class of order $p^n \geq p^{p+1}$ and exponent $p^e$, where $p$ is an odd prime. Denote by $\Omega_0(p) := {\mathbb N}$ and for $N \in {\mathbb N}$, let $\Omega_N(p)$ is the set of all solutions $y$ of the linear diophantine equation 
\begin{equation*}
y = hp^{N} + \sum_{i=1}^{N} {\frac {1}{2}} (p^N - p^{N-i})x_i
\end{equation*}
in integer variables $h, x_1, \dotsc, x_{N} \geq 0$. Then the genus spectrums of $G$ given by 
\[
{\mathrm{sp}}(G) = p^{n-e} {\widetilde{\mathrm{sp}}}(G) + 1,
\hspace*{.5in} {\mathrm{where}}~~~ {\widetilde{\mathrm{sp}}}(G) = \bigcup_{N=0}^{e} p^{e-N} F_{N} + p^{e-N} \Omega_N(p)
\]
where each $F_i$ is a finite set of integers listed in the tables (\ref{Fi-sets-table1}), (\ref{Fi-sets-table2}) and (\ref{Fi-sets-table3}).

\smallskip

\noindent From table (\ref{Fi-sets-table1}), we see that while $|G| \geq p^{2p}$, there are only five genus spectra, as it follows from {\ref{exp-p-p2-lemma}}, that exactly one of the conditions ${\mathrm{ed}}(G) \geq 2$, or ${\mathrm{ed}}(G) = 1$ hold for a fixed $|G|$. Table (\ref{Fi-sets-table2}) determines the spectra of groups with $p^{p+2} \leq |G| \leq p^{2p-1}$, which satisfy ${\mathrm{ed}}(G) \geq 2$ by {\ref{exp-p-p2-lemma}}. Here we have seven spectra with one generic case: Type V, the only group of order $243$, namely ${\mathrm{SmallGroup}}(243,28)$, checked by GAP {\cite{gap}}. Table (\ref{Fi-sets-table3}) shows that there are eight distinct spectra while $|G| = p^{p+1}$, with a generic case: Type IV, the only group of order $81$, which is the wreath product $C_3 \wr C_3$ (also ${\mathrm{SmallGroup}}(81,7)$), with its two step centralizer of exponent $3$. 
 
\smallskip

\noindent A combinatorial description of the sets $\Omega_N(p)$ is given in \cite{kma}. A reduced stable upper genus of a finite group $G$ is a unique (scaled by the factor $p^{n-e}$) least integer so that all integers beyond this is a reduced genus (see \ref{basic-action}). Using these descriptions of the signatures we proved that :

\smallskip

\noindent {\bf Theorem B.} The stable upper genus of a finite $p$-group of maximal class of order $\geq p^{p+2}$ and exponent $p^e \geq p^3$ is a constant $\sigma_e(p)-p^e$ (i.e., independent of $|G|$) if and only if 

\smallskip

(i) $p \geq 5$ while ${\mathrm{ed}}(G) \geq 2$, and
(ii) $p \geq 7$ while ${\mathrm{ed}}(G) = 1$.   

\smallskip

\noindent The paper is organized as follows : Section {\ref{prelim}} discuss all the basic results and some useful inequalities connecting the order and exponent of the $p$-groups of maximal classes. Section {\ref{z-classes}} discuss the notion of $z$-classes and the relevant results relating the notion of $p$-group of maximal class of exceptional types. Section {\ref{resexp3ed2}} and {\ref{resed1}} describe the group actions in the cases ${\mathrm{ed}}(G) \geq 2$ and ${\mathrm{ed}}(G)=1$ respectively while $|G| \geq p^{p+2}$. Section {\ref{orderp^p+1}} describe the actions while $|G| = p^{p+1}$. Using these results, in section {\ref{spectrums}} we compute the minimum genus and classify all groups which has stable upper genus $\sigma_e(p)-p^e$. The final section {\ref{examples}} construct the examples of the $p$-groups of maximal class that achieve the spectrums.


\begin{table}[H]
  \begin{center}
  \caption {$e \geq 3, |G| \geq p^{2p}$}\label{Fi-sets-table1}  
 \begin{adjustbox}{angle=90}
 \small
\begin{tabular}{| m{2em} | m{2.5em} | m{18cm} |}
\hline
 Type & ${\mathrm{ed}}(G)$ & $F_i$  \\ \hline
 I    & $\geq 2$ & $F_e = \left\{ p^e, {\frac {p^e - 1}{2}}, {\frac {p^e - 2p^{e-1} - 1}{2}} \right\}, F_1 = \left\{ p-1, {\frac {p-3}{2}} \right\}, F_0 = \emptyset$, \newline $F_N = \left\{ {\frac {ap^N - bp^{N-1} - 1}{2}} : (a,b) = (3,2), (2,3) \right\}~ (2 \leq N \leq e-1)$  \\ \cline{2-3}
 	  & $1$ & $F_e = \left\{ p^e - 1, {\frac {3p^e - 2p^{e-1} - 1}{2}}, {\frac {p^e - 2p^{e-1} - 1}{2}} \right\}, F_1 = \left\{ p - 1, {\frac {3p - 1}{2}}, {\frac {p - 3}{2}} \right\}, F_0 = \{ 0 \}$, \newline $F_N = \left\{ {\frac {3p^N - 1}{2}}, {\frac {3p^N - 2p^{N-1} - 1}{2}}, {\frac {2p^N - 3p^{N-1} - 1}{2}} \right\}~ (2 \leq N \leq e-1)$ \\ \hline 
 II   & $\geq 2$ & $F_e = \left\{ p^e, {\frac {p^e - 1}{2}}, {\frac {p^e - 2p^{e-2} - 1}{2}} \right\}, F_2 = \left\{ p^2 - 1, {\frac {p^2 - 3}{2}} \right\}, F_1 = \emptyset = F_0$, \newline $F_N = \left\{ {\frac {ap^N - bp^{N-2} - 1}{2}} ~:~ (a,b) = (3,2), (2,3) \right\}~ (2 \leq N \leq e-1)$ \\ \cline{2-3}
 	  & $1$ & $F_e = \left\{ p^e - 1, {\frac {3p^e - 2p^{e-2} - 1}{2}}, {\frac {p^e - 2p^{e-2} - 1}{2}} \right\}, F_2 = \left\{ p^2 - 1, {\frac {3p^2 - 1}{2}}, {\frac {p^2 - 3}{2}} \right\}, F_1 = \left\{ {\frac {3p - 1}{2}} \right\}, F_0 = \{ 0 \}$,
$F_N = \left\{ {\frac {3p^N - 1}{2}}, {\frac {3p^N - 2p^{N-2} - 1}{2}}, {\frac {2p^N - 3p^{N-2} - 1}{2}} \right\}~ (3 \leq N \leq e-1)$ \\ \hline
 III  & $\geq 2$ & $F_e = \left\{ p^e, {\frac {p^e - 1}{2}}, {\frac {p^e - ap^{e-1} - bp^{e-2} - 1}{2}} : (a,b) = (0,2), (1,1), (2,0) \right\}$, \newline $F_N = \left\{ {\frac {ap^N - bp^{N-1} - cp^{N-2} - 1}{2}} : a = 2,3, a+b+c = 5, b,c \geq 0, (b,c) \neq (1,1) \right\} (3 \leq N \leq e-1)$, $F_2 = \left\{ p^2 - 1, {\frac {3p^2 - 2p - 1}{2}}, {\frac {ap^2 - bp - c}{2}} ~:~ b+c=3, b \geq 0, c \geq 1 \right\}, F_1 = \left\{ {\frac {p-3}{2}} \right\}, F_0 = \emptyset$ \\ \cline{2-3}
 \parbox[t]{2mm}{\multirow{3}{*}{\rotatebox[origin=c]{270}{$p \geq 5$}}} & $1$ & $F_e = \left\{ p^e - 1, {\frac {ap^e - bp^{e-1} - cp^{e-2} - 1}{2}} : a = 1,3, (b,c) = (2,0), (1,1), (0,2) \right\}$, \newline $F_N = \left\{ {\frac {ap^N - bp^{N-1} - cp^{N-2} - 1}{2}} : (a,b,c) = (3,0,0), (3,2,0),(3,0,2), (2,3,0), (2,2,1), (2,1,2), (2,0,3) \right\}$, \newline $(3 \leq N \leq e-1)$, \newline $F_2 = \left\{ p^2 - 1, {\frac {3p^2 - 1}{2}}, {\frac {3p^2 - 2p - 1}{2}}, {\frac {p^2 - bp - c}{2}} ~:~ (b,c) = (2,1), (1,2), (0,3) \right\}, F_1 = \left\{ p - 1, {\frac {3p - 1}{2}}, {\frac {p-3}{2}} \right\}, F_0 = \{ 0 \}$ \\ \hline
 IV   &  $\geq 2$ & $F_e = \left\{ p^e, {\frac {p^e - 1}{2}}, {\frac {p^e - p^{e-1} - p^{e-2} - 1}{2}}, {\frac {p^e - 2p^{e-2} - 1}{2}} \right\}$, $F_2 = \left\{ p^2 - 1, {\frac {3p^2 - 2p - 1}{2}}, {\frac {p^2 - p - 2}{2}} \right\}, F_1 = \left\{ p - 1 \right\}, F_0 = \emptyset$, \newline $F_N = \left\{ {\frac {ap^N - bp^{N-1} - cp^{N-2} - 1}{2}} : (a,b,c) = (3,0,2), (3,2,0), (2,1,2) \right\}, (3 \leq N \leq e-1)$, \\ \cline{2-3}
 	  & $1$ & $F_e = \left\{ p^e - 1, {\frac {ap^e - bp^{e-1} - cp^{e-2} - 1}{2}} : a = 1,3, (b,c) = (1,1), (0,2) \right\}, F_2 = \left\{ p^2 - 1, {\frac {3p^2 - 1}{2}}, {\frac {3p^2 - 2p - 1}{2}}, {\frac {p^2 - p - 2}{2}} \right\}$, \newline $F_N = \left\{ {\frac {ap^N - bp^{N-1} - cp^{N-2} - 1}{2}} : (a,b,c) = (3,0,0), (3,2,0),(3,0,2), (2,1,2) \right\}, (3 \leq N \leq e-1)$, \newline $F_1 = \left\{ p - 1, {\frac {3p - 1}{2}} \right\}, F_0 = \{ 0 \}$ \\ \hline
 V    &  $\geq 2$ & $F_e = \left\{ p^e, {\frac {p^e - 1}{2}}, {\frac {p^e - p^{e-1} - p^{e-2} - 1}{2}}, {\frac {p^e - 2p^{e-1} - 1}{2}} \right\}$, $F_2 = \left\{ p^2 - 1, {\frac {3p^2 - 2p - 1}{2}}, {\frac {p^2 - 2p - 1}{2}} \right\}, F_1 = \left\{ p - 1 \right\}, F_0 = \emptyset$, \newline $F_N = \left\{ {\frac {ap^N - bp^{N-1} - cp^{N-2} - 1}{2}} : (a,b,c) = (3,0,2), (3,2,0), (2,2,1) \right\}, (3 \leq N \leq e-1)$ \\ \cline{2-3}
	 & $1$ &  $F_e = \left\{ p^e - 1, {\frac {ap^e - bp^{e-1} - cp^{e-2} - 1}{2}} : a = 1, 3, (b,c) = (1,1), (2,0) \right\}$, \newline $F_N = \left\{ {\frac {ap^N - bp^{N-1} - cp^{N-2} - 1}{2}} : (a,b,c) = (3,0,0), (3,2,0), (3,0,2), (2,2,1) \right\}, (3 \leq N \leq e-1)$, \newline $F_2 = \left\{ p^2 - 1, {\frac {3p^2 - 1}{2}}, {\frac {3p^2 - 2p - 1}{2}}, {\frac {p^2 - 2p - 1}{2}} \right\}, F_1 = \left\{ p - 1, {\frac {3p - 1}{2}} \right\}, F_0 = \{ 0 \}$ \\ \hline
\end{tabular}
\end{adjustbox}
\end{center}
\end{table}

\begin{table}[H]
  \begin{center}
  \caption {$e = 2$, ${\mathrm{ed}}(G) \geq 2$, $p^{p+2} \leq |G| \leq p^{2p-1}$}\label{Fi-sets-table2}  
\begin{tabular}{| m{0.9cm} | m{2.8cm} | m{13.6cm} |}
\hline
 Type & Property &  $F_i$  \\ \hline
 I    & -        & $F_2 = \left\{ p^2, {\frac {p^2 - 1}{2}}, {\frac {p^2 - 2p - 1}{2}} \right\}, F_1 = \left\{ p-1, {\frac {p-3}{2}} \right\}, F_0 = \emptyset$ \\ \hline
 II   & -        & $F_2 = \left\{ p^2, {\frac {p^2 - 1}{2}}, {\frac {p^2 - 3}{2}} \right\}, F_1 = \emptyset = F_0$ \\ \hline
 III  & $(p,p,p)$-type & $F_2 = \left\{ p^2, {\frac {p^2 - 1}{2}}, {\frac {p^2 - ap - b}{2}} : (a,b) = (2,1), (1,2), (0,3) \right\}$, $F_1 = \left\{ p-1, {\frac {p-3}{2}} \right\}$, $F_0 = \emptyset$ \\ \cline{2-3}
  $p \geq 5$ & not $(p,p,p)$-type & $F_2 = \left\{ p^2, {\frac {p^2 - 1}{2}}, {\frac {p^2 - ap - b}{2}} : (a,b) = (2,1), (1,2), (0,3) \right\}, F_0 = \emptyset$,  $F_1 = \left\{ p-1, p-2 \right\}$ \\ \hline
  IV  &  -       & $F_2 = \left\{ p^2, {\frac {p^2 - 1}{2}}, {\frac {p^2 - ap - b}{2}} : (a,b) = (1,2), (0,3) \right\}, F_1 = \left\{ p-1 \right\}, F_0 = \emptyset$ \\ \hline
  V   &  $p \geq 5$   & $F_2 = \left\{ p^2, {\frac {p^2 - 1}{2}}, {\frac {p^2 - 2p -1}{2}} \right\}, F_0 = \emptyset, F_1 = \left\{ p-1, {\frac {p-3}{2}} \right\}$ \\ \cline{2-3}
      &  $p = 3$   & $F_2 = \left\{ 9, 4, 1 \right\}, F_0 = \emptyset, F_1 = \left\{ 2, 1 \right\}$ \\ \hline
\end{tabular}
\end{center}
\end{table}

\noindent Tables {\ref{5^6-grps}} and {\ref{5^7-grps}} provide a list of groups ${\mathrm{SmallGroup}}(5^n,i)$ of order $5^n, (n = 6,7)$ along with their types which are computed using GAP {\cite{gap}}. Based on the observations we can ask several questions :

\smallskip

\noindent {\bf Co-class wise uniformity conjecture :} Let ${\mathcal G}_{n,e,c}$ denote the set of all finite $p$-groups $G$ of co-class $c := n-cl(G)$ ($cl(G)$ denote the nilpotency class of $G$) with order $|G| = p^n$ and exponent $p^e$. Does there exists a polynomial function $f(c)$, depending only on $c$, so that the number of genus spectrums of $G \in {\mathcal G}_{n,e,c}$ is bounded above by $f(c)$? 

\noindent The results of this paper shows $f(1) \leq 8$. If we distinguish the the primes, we proved $f(1,3) \leq 8$ and $f(1,p) \leq 7$ for any prime $p \geq 5$ (while $f(1,2) = 3$ was already known from unpublished thesis of Yasemin Talu, 1993).

\begin{table}[H]
  \begin{center}
  \caption {$e = 2$, $|G| = p^{p+1}$}\label{Fi-sets-table3}
  \small  
\begin{tabular}{| m{0.72cm} | m{1.2cm} | m{2.6cm} | m{12.15cm} |}
\hline
 Type & ${\mathrm{exp}}(G_1)$ & Property &  $F_i$  \\ \hline
  I   & $p^2$ & - & $F_2 = \left\{ p^2 - 1, {\frac {3p^2 - 2p - 1}{2}}, {\frac {p^2 - 2p - 1}{2}} \right\}, F_1 = \left\{ {\frac {p-1}{2}}, {\frac {p-3}{2}} \right\}, F_0 = \{ 0 \}$ \\ \hline
  II  &  $p$ & - & $F_2 = \left\{ p^2 - 1, {\frac {p^2 - ap - b}{2}} ~:~ (a,b) = (0,3), (1,2) \right\}, F_1 = \left\{ {\frac {p-1}{2}} \right\}, F_0 = \{ 0 \}$ \\ \cline{2-4}
      &  $p^2$ & - & $F_2 = \left\{ p^2 - 1, {\frac {p^2 - 3}{2}} \right\}, F_1 = \left\{ {\frac {p-1}{2}} \right\}, F_0 = \{ 0 \}$ \\ \hline
  III & $p$ & $(p^2,p^2,p^2)$-type & $F_2 = \left\{ p^2 - 1, p^2 - 2, {\frac {3p^2 - 2p - 1}{2}}, {\frac {p^2 - ap - b}{2}} ~:~ (a,b) = (0,3), (1,2),(2,1) \right\}$, \\
      &     &    &  $F_1 = \left\{ {\frac {p-1}{2}}, {\frac {p-3}{2}} \right\}, F_0 = \{ 0 \} $ \\ \cline{3-4}
   \parbox[t]{2mm}{\multirow{3}{*}{\rotatebox[origin=c]{90}{$p \geq 5$}}}   &     & not & $F_2 = \left\{ p^2 - 1, p^2 - 2, {\frac {3p^2 - 2p - 1}{2}}, {\frac {p^2 - ap - b}{2}} ~:~ (a,b) = (1,2),(2,1) \right\}$, \\ 
   &  & $(p^2,p^2,p^2)$-type &  $F_1 = \left\{ {\frac {p-1}{2}}, {\frac {p-3}{2}} \right\}, F_0 = \{ 0 \} $ \\ \cline{2-4}  
      & $p^2$ & $(p,p,p)$-type & $F_2 = \left\{ p^2 - 1, {\frac {3p^2 - 2p - 1}{2}}, {\frac {p^2 - ap - b}{2}} ~:~ (a,b) = (0,3), (1,2),(2,1) \right\},$ \\ 
      &       &                & $F_1 = \{ {\frac {p-1}{2}}, {\frac {p-3}{2}} \}, F_0 = \{ 0 \}$ \\ \cline{3-4}
      &       & not            & $F_2 = \left\{ p^2 - 1, {\frac {3p^2 - 2p - 1}{2}}, {\frac {p^2 - ap - b}{2}} ~:~ (a,b) = (0,3), (1,2),(2,1) \right\},$ \\
      &       & $(p,p,p)$-type & $F_1 = \{ p-2, {\frac {p-1}{2}} \},  F_0 = \{ 0 \}$ \\ \cline{1-3}
  IV  &	$p$  & $p \geq 5$      &	  \\ \cline{3-4}
      &      & $p = 3$         & $F_2 = \left\{ 1, 2, 7, 8, 10 \right\}, F_1 = \left\{ 1 \right\}, F_0 = \{ 0 \}$ \\ \cline{2-4} 
      & $p^2$ &    & Same as Type II, ${\mathrm{exp}}(G_1) = p$ \\ \hline 
 V	  &  $p$ &      &  Same as Type I \\ \cline{2-4}
 	  &  $p^2$ & $p \geq 5$ & Same as Type III, ${\mathrm{exp}}(G_1) = p$, not $(p^2,p^2,p^2)$-type \\ \cline{3-4}
 	  &      &  $p = 3$ & Same as Type IV, ${\mathrm{exp}}(G_1) = 3$ \\ \hline       
\end{tabular}
\end{center}
\end{table}


\begin{table}[H]
  \begin{center}
  \caption {$e = 2$, ${\mathrm{ed}}(G) \geq 2$, $|G| = 78125$}\label{5^7-grps}  
  \small
\begin{tabular}{| m{0.9cm} | m{2.8cm} | m{9.6cm} | m{2.6cm}|}
\hline
 Type & Property &  $i$  & No. of groups \\ \hline
 I    & -        &  $1283, 1286, 1297, 1304, 1370$ & $5$ \\ \hline
 II   & -        &  $1284, 1287-1290, 1292, 1298-1301, 1305, 1308, 1310, 1312$, & \\ \cline{3-3} 
 	  &          &  $1313, 1318, 1322, 1323,1327,1328, 1330,1331, 1337,1338$, & $39$ \\ \cline{3-3} 
 	  &          &  $1342,1343, 1345, 1346, 1352, 1353, 1355, 1358,1361,1363$, & \\ \cline{3-3} 
 	  &          &  $1371,1372, 1375, 1377, 1379$ & \\ \hline
 III  & $(5,5,5)$-type & - & $0$ \\ \cline{2-4}
   & not $(5,5,5)$-type &  $1291,1306,1307,1311,1314,1316,1319,1321,1324,1326$ & \\ \cline{2-3}
      &          & $1329,1332,1334,1336,1339,1341,1344,1347,1349,1351$ & $27$ \\ \cline{2-3}
      &          & $1354,1357,1360,1364,1373,1374,1378$ & \\ \hline
  IV  &  -       & $1282,1285,1293-1296,1302,1303,1309,1315,1317,1320$ & \\ \cline{3-3} 
 	  &          & $1325,1333,1335,1340,1348,1350,1356,1359,1362,$ & $28$ \\ \cline{3-3} 
 	  &          &  $1365-1369,1376,1380$ & \\ \hline
  V   &  -      & - & $0$ \\ \hline
\end{tabular}
\end{center}
\end{table}
 
\smallskip

\begin{table}[H]
  \begin{center}
  \caption {$e = 2$, $|G| = 15625$}\label{5^6-grps}  
  \small
\begin{tabular}{| m{0.72cm} | m{1.2cm} | m{3.4cm} | m{7.1cm} | m{2.5cm} |}
\hline
 Type & ${\mathrm{exp}}(G_1)$ & Property &  $i$ & No. of groups \\ \hline
  I   & $25$ & - & $630,636,652$ & $3$ \\ \hline
  II  &  $5$ & - &  $644,662$ & $2$ \\ \cline{2-5}
      &  $25$ & - & $635,637,638,648,650,660,666,668$ & $8$ \\ \hline
  III & $5$ & $(25,25,25)$-type &  $643,645,661,663$ & $4$ \\ \cline{3-5}
      &     & not $(25,25,25)$-type & - & $0$ \\ \cline{2-5}  
      & $25$ & $(5,5,5)$-type & - & $0$ \\ \cline{3-5}
      &       & not $(5,5,5)$-type & $646,647,649,664,665,667$ & $6$ \\ \hline
  IV  &	$5$  & -      &	 $631,639,651,656$ & $4$ \\ \cline{2-5} 
      & $25$ &    & $632-634,640-642,653-655,657-659$ & $12$ \\ \hline 
 V	  &  $5$ &      &  - & $0$ \\ \cline{2-5}
 	  &  $25$ & - & - & $0$ \\ \hline       
\end{tabular}
\end{center}
\end{table}

\noindent {\bf Existence type questions :} are there groups of (i) $p^2$-exceptional types? (ii) $(p,p,p)$-type in case $G$ is not exceptional. 

\smallskip

\noindent Thm.{\ref{metab-non-p2}} and \cite[Thm.2.2]{fggs} shows that one need to understand the co-class $1$ groups with $cl(G_1) \geq 3$, which is in general hard (see {\cite{mie2}}). 

\bigskip

\noindent {\bf Notations}

$\bullet$ ${\mathbb N}_{\geq k} = \left\{ n \in {\mathbb N} ~;~ n \geq k \right\}$,  

$\bullet$ For $\alpha, \beta \in {\mathbb R}; A, B \subseteq {\mathbb R}$; $\alpha A + \beta B := \{ \alpha a + \beta b ~:~ a \in A, b \in B \}$ \\
\hspace*{3in} and is $\emptyset$ if one of $A, B$ is empty, 

$\bullet$  $C_n$ denote the cyclic group of order $n$,

$\bullet$ ${\mathrm{exp}}(G) := {\mathrm{l.c.m.}} \{ |x| ~:~ x \in G \}$ is exponent of a group $G$,

$\bullet$ For a finite $p$-group $G$ of maximal class, ${\mathrm{ed}}(G)$ is its exponent depth (Definition \ref{edG-def}),

$\bullet$  $[H, K] = \langle [x, y] = x^{-1} y^{-1} xy ~:~ x \in H, y \in K \rangle$, where $H, K \leq G$, 

$\bullet$  $\gamma_1(G) = G_0 = G, G_{i+1} = \gamma_{i+1}(G) = [\gamma_i(G), G]$ for $i \geq 1$, $G_1 = C_{G}(G_{2}/G_{4})$,

$\bullet$ $G^{\ast}_i = G_i \setminus G_{i+1} \hspace{.2in} (1 \leq i \leq n-1)$,

$\bullet$ $Z(G)$ is the center of $G$,

$\bullet$  $cl(G)$ is the nilpotency class of $G$, 

$\bullet$ For $X \subset G$, $||X|| = \{ |x| ~:~ x \in X \}$,

$\bullet$ For any $0 \leq i \leq n-1$, $G^{\ast}_i = G_i \setminus G_{i+1}$.

$\bullet$  For a finite $p$-group, $\mho_1(G) = \langle x^p ~:~ x \in G \rangle$. 

$\bullet$  For a finite $p$-group of maximal class, $c(G)$ denote the degree of commutativity of $G$.  

$\bullet$ $H(\Lambda)$ (resp. $E_i(\Lambda)$) are the set of hyperbolic (resp. elliptic generators of order $p^i$) from $\Lambda$. \\
\hspace*{.25in} $E(\Lambda) = \cup_{i=1}^{e} E_i(\Lambda)$. 

$\bullet$ ${\mathrm{ed}}(G)$ denote the exponent depth of $G$ (Definition \ref{edG-def}).

\section{\bf Preliminaries}
\label{prelim}

\subsection{Groups acting on surfaces (A2.1)}\label{basic-action} 

The results of this section are true for arbitrary finite groups $G$ (see {\cite{kul}}). We restrict ourself to the case while $G$ is a group with order $p^n$ and exponent $p^e$, where $p$ is an odd prime, $n \geq e \geq 1$ are integers. Then $G$ acts on $\Sigma_g$ with orbit space $\Sigma_h \cong \Sigma_g/G$ if and only if 
\begin{equation}\label{ec-eqn}
2(g-1) = p^n \Bigl\{ 2(h-1) + \sum_{i=1}^{e} m_i \bigl( 1 - {\frac {1}{p^i}} \bigr) \Bigr\}
\end{equation}
and there exists a subset $\Lambda = \{ a_1, b_1, \dotsc, a_h, b_h, X_{ij} ~(1 \leq i \leq e, 1 \leq j \leq m_i) \} \subseteq G$ with $G = \langle \Lambda \rangle$ satisfying conditions
\begin{equation}\label{fuch-gen}
|X_{ij}| = p^i, \hspace{.2in} \prod_{k=1}^{h} [a_k, b_k] \prod_{i=1}^{e} \prod_{j=1}^{m_i} X_{ij} = 1 \hspace{1em} {\mathrm{(Long~relation)}}
\end{equation}
If such happens we call the $(e+1)$-tuple $\sigma = (h; m_1, \dotsc, m_e)$ a signature of $G$. In the above partial presentation of $G$, we denote 
\[
H(\Lambda) = \{ a_1, b_1, \dotsc, a_h, b_h \}, ~ E_i(\Lambda) = \{ X_{ij} ~:~ 1 \leq j \leq m_i \} ~(1 \leq i \leq e), ~ E(\Lambda) = \bigcup_{i=1}^{e} E_i(\Lambda)
\]
and the elements of $H(\Lambda)$ (resp. $E_i(\Lambda)$) are called the {\bf hyperbolic generators} (resp. {\bf elliptic generators}). The relation in ({\ref{fuch-gen}}) above involving the product is called the {\bf long relation}. Let ${\mathcal{S}}(G)$ denote the collection of all such signatures of $G$, referred as the {\bf signature spectrum} of $G$. 

If $\sigma = (h; m_1, \dotsc, m_e) \in {\mathcal{S}}(G)$ we denote by $N = N(\sigma)$ as
\[
N(\sigma) = 
\begin{cases}
0 &\mbox{if } m_i = 0 {\mathrm{~for~all~}} 1 \leq i \leq e \\
k &\mbox{if } m_k \neq 0, m_i = 0 {\mathrm{~for~all~}} i > k
\end{cases}
\] 

Now consider the genus map ${\mathrm{\bf g}} : {\mathcal{S}}(G) \rightarrow {\mathrm {sp}}(G)$ defined by
\[
{\mathrm{\bf {g}}}(h; m_1, \dotsc, m_e) = p^{n-e} \bigl\{ (h-1)p^e + \sum_{i=1}^{e} {\frac {1}{2}}(p^e - p^{e-i}) m_i \bigr\} + 1
\]
The scaled map 
\[
{\tilde{\mathrm{\bf {g}}}}(h; m_1, \dotsc, m_e) = (h-1)p^e + \sum_{i=1}^{e} {\frac {1}{2}}(p^e - p^{e-i}) m_i
\]
is called the reduced genus map. While $\sigma \in {\mathcal S}(G)$ and ${\mathrm{\bf {g}}}(\sigma)$ is the corresponding genus, we refer ${\tilde{\mathrm{\bf {g}}}}(\sigma)$ as the corresponding reduced genus. The set of all reduced genus is denoted by ${\widetilde{\mathrm{sp}}}(G)$. Following the results of {\cite{kul}} we have that while $p$ is an odd prime, the image of the reduced genus map ${\tilde{\mathrm{\bf {g}}}}({\mathcal{S}}(G)) = {\widetilde{\mathrm{sp}}}(G)$ is a co-finite subset of $\N$. Similar result holds for $p=2$, which require modification of the factor $2^{n-e}$ that depends on certain special property of $G$. The number $\mu(G) := \min{\mathrm{sp}}(G)$ is called the {\bf minimum genus} of $G$ and its scaled factor denoted as ${\tilde{\mu}}(G)$, is called the reduced minimum genus. Since $|\N \setminus {\widetilde{\mathrm{sp}}}(G)| < \infty$, it contain ${\tilde{\sigma}}(G)$, called the {\bf reduced stable upper genus} which is minimum among all ${\tilde{g}} \in {\widetilde{\mathrm{sp}}}(G)$ so that ${\N}_{\geq {\tilde{g}}} \subseteq {\widetilde{\mathrm{sp}}}(G)$. Scaling back, we have $\sigma(G) = p^{n-e} {\tilde{\sigma}}(G) + 1$, which is called the stable upper genus of $G$.

\subsection{Some basic principles (A2.2)}\label{basic} Here we discuss some basic principles which are followed for generation of finite group actions on surfaces. 

\noindent (i) {\it conjugacy principle :} Suppose the group $G$ acts on $\Sigma_g$ satisfying the equations (\ref{ec-eqn}) and (\ref{fuch-gen}) as above except the elliptic generators in the long relation are misplaced. Then we can modify these generators using conjugacy. For example if $Y_1 Y_2$ appears in the product (as one of the $X_{ij}$'s). Then replace the generators $Y_1, Y_2$ by $Y_2^{Y^{-1}_1}, Y_1$ and the product by $(Y_1 Y_2 Y_1^{-1}) Y_1$. These doesn't change the long relation in $G$ and the group is still generated by the new set of generators. Hence in finitely many steps one can obtain the long relation in the correct form as in (\ref{fuch-gen}).

\noindent (ii) {\it extension principle :} Suppose $p$ an odd prime and $\sigma = (h; m_1, \dotsc, m_N, 0, \dotsc, 0) \in {\mathcal S}(G)$ with $m_N \geq 1$ realized by generators satisfying (\ref{fuch-gen}). Let $Y$ denote an elliptic generator of order $p^N$. Then for any $1 \leq k \leq N$, $Y$ can be written as $Y = Y_1 Y_2$ with $|Y_1| = p^k, |Y_2| = p^N$. Replacing the generator $Y$ by $Y_1, Y_2$ the group $G$ is still generated and we can simply replace $Y$ by $Y_1 Y_2$ in the long relation. Finally using conjugacy principle we have another signature which differ from $\sigma$ by changing $m_k$ to $m_k + 1$. Moreover for any $x \in G$ we can add a pair of generators both equal to $x$ and insert the trivial product $[x, x]$ in the long relation. Thus any $\sigma^{\prime} = (h^{\prime}; m_1^{\prime}, \dotsc, m_N^{\prime}, 0, \dotsc, 0) \in {\mathcal S}(G)$ with $h^{\prime} \geq h$ and $m_j^{\prime} \geq m_j$ for each $j$. 

\subsection{Finite $p$-groups of maximal class} 

A finite $p$-group $G$ of order $p^n, (n \geq 4)$ is said to be of maximal class (or co-class $1$) if its nilpotency class $cl(G)$ is $n-1$. In this case we have 
\[
G/{\gamma_2(G)} \cong C_p \times C_p, \gamma_i(G)/{\gamma_{i+1}(G)} \cong C_p ~~{\mathrm {for}}~~ 2 \leq i \leq n - 1
\]
The two step centralizers $C_G(\gamma_i(G)/{\gamma_{i+2}(G)})$ induced by the natural conjugacy action of $G$ on the quotient group $\gamma_i(G)/{\gamma_{i+2}(G)}$ are maximal subgroups of $G$. If we denote $G_0 := G, G_1 = C_G(\gamma_2(G)/{\gamma_{4}(G)})$ and reset the definitions $G_i := \gamma_i(G)$ for $i \geq 2$, then $G_i/G_{i+1} \cong C_p$ for each $0 \leq i \leq n-1$. Most of the following results are standard and we refer {\cite{falc}} and {\cite{lmc}} for the proof of these. 

If $G$ is $p$-group of maximal class, it is well known that it has at the most two two-step centralizers, namely $G_1$ and $C_G(G_{n-2})$. If the degree of commutativity $c(G)=0$ we have $G_1 \neq C_G(G_{n-2})$ and they coincide while $c(G) \geq 1$. For small order we have the standard result :

\begin{lemma}\cite[Thm. 3.3.2]{lmc}
\label{p+1-lemma}  Let $G$ be a $p$-group of maximal class of order $\leq p^{p+1}$. Then ${\mathrm{exp}}(G/Z(G)) = {\mathrm{exp}}(G_2) = p$. 
\end{lemma}

\noindent The elements $s \in G \setminus (G_1 \cup C_G(G_{n-2}))$ are called the {\bf uniform elements}. From Blackburn's theorem (see \cite{bla}) it follows that $G$ always contain uniform elements. The following standard result will be useful for us.

\begin{proposition}\cite[Thm. 3.15]{falc}
\label{unif-prop}  Let $G$ be a $p$-group of maximal class of order $p^n$ and $s \in G \setminus (G_1 \cup C_G(G_{n-2}))$. Then we have 
\begin{enumerate}[(i)]
\item $C_G(s) = \langle s \rangle Z(G)$.

\item $s^p \in Z(G)$ and hence $|s| \leq p^2$. Consequently $|C_G(s)| = p^2$.

\item The conjugacy class $s^G$ is precisely equal to the coset $s G_2$. 
\end{enumerate}
\end{proposition}

\noindent A $p$-groups $G$ of maximal class is minimally generated by two elements, $G = \langle s, s_1 \rangle$, where $s$ is an uniform element and $s_1 \in G_1 \setminus G_2$. Set $s_0 := s$ and $s_{i+1} := [s_i, s]$ for $i \geq 1$, so that $s_i \in G_i \setminus G_{i+1}$. The finite $p$-groups of maximal class of order $\leq p^p$ are known to be regular, while if its order is $\geq p^{p+1}$, they are necessarily irregular and hence has exponent $\geq p^2$. We summarize some of these properties in the following. For proof of these see \cite[Thm. 4.9, Exercise 4.1]{falc}.

\begin{theorem}
\label{powstruc} Let $G$ be a $p$-group of maximal class of order $p^n \geq p^{p+1}$. Then we have :
\begin{enumerate}[(i)]
\item $G$ is irregular with regular maximal subgroup $G_1$,

\item While $|G| \geq p^{p+2}$, $\mho_1(G_i) = G_{i+p-1}$ for every $i \geq 1$, 

\item For each $1 \leq i \leq n-p$ and $x \in G_i \setminus G_{i+1}$ we have $x^p \in G_{i+p-1} \setminus G_{i+p}$. 
\end{enumerate}
\end{theorem}

\noindent A consequence of the above theorem using the commutator collection formula is the following:

\begin{proposition} \cite[Prop. 3.3.8]{lmc} If $c(G) \geq 1$, then $s^p_i s_{i+p-1} \in G_{i+p}$ for each $i \geq 2$ and $s^p_1 s_p \in G_{p+1} G_{n-1}$.  
\end{proposition}  

\begin{lemma} \cite[Exercise 2, \textsection 9]{ber} Let $G$ be a $p$-group of maximal class of order $p^n \geq p^{p+2}$. Then ${\mathrm{exp}}(G) = {\mathrm{exp}}(G_1)$.
\end{lemma}

The following result connect the order $p^n$ and the exponent $p^e$ of $G$ by some useful inequality. 

\begin{theorem}\label{exp-p-p2-lemma}
Let $G$ be a $p$-group of maximal class of order $p^n \geq p^{p+2}$ and exponent $p^e$. Then we have :

\begin{enumerate}[(1)]

\item $||G^{\ast}_{n-p}|| = \{ p^2 \}$ and ${\mathrm{exp}}(G_{n-p+1}) = p$, 

\item For any $1 \leq i \leq n-1$, let $k(i)$ be the unique integer $\geq -1$ satisfying 
\[
i + k(i)(p-1) \leq n-p < i + (k(i)+1)(p-1)
\]
then $||G^{\ast}_i|| = \{ p^{k(i)+2} \}$.

\item $e \geq 2$ and \[
1+(e-2)(p-1) \leq n-p < 1+(e-1)(p-1)
\] 
Conversely, if $G$ is a $p$-group of maximal class of order $p^n \geq p^{p+2}$, where $n$ satisfies the above inequality for some positive integer $e \geq 2$, then ${\mathrm{exp}}(G) = p^e$. When this happens, $||G^{\ast}_1|| = \{ p^e \}$.

\item 
\[
||G^{\ast}_i|| = \{ p^N \} ~~{\mathrm{for~each}}~~ n - N(p-1) \leq i \leq n - (N-1)(p-1) - 1
\]
for each $N \leq e-1$, and 
\[
||G^{\ast}_i|| = \{ p^e \} ~~{\mathrm{for~each}}~~ 1 \leq i \leq n - (e-1)(p-1) - 1
\]

\end{enumerate}

\end{theorem}

\noindent {\bf Proof.} (1) If $x \in G^{\ast}_{n-p}$ we have by {\ref{powstruc}}(iii) $x^p \in G^{\ast}_{n-1}$. Thus $x^p \neq 1$ and $x^{p^2} = 1$. Next by {\ref{powstruc}}(ii) $\mho_1(G_{n-p+1}) = 1$. 

(2) For $x \in G^{\ast}_i$, repeated use of {\ref{powstruc}}(iii) yield $x^{p^{k(i)}} \in G^{\ast}_{i+k(i)(p-1)}$ and $x^{p^{k(i)+1}} \in G^{\ast}_{i+(k(i)+1)(p-1)}$. Thus $x^{p^{k(i)+1}}$ is a non-trivial element of $G_{i+(k(i)+1)(p-1)} \subseteq G_{n-p+1}$. Thus by {\ref{exp-p-p2-lemma}}, $|x| = p^{k(i)+2}$. 

(3) Since exp($G$) $=$ exp($G_1$) and $G_1$ is regular, $G^{\ast}_1$ must contain an element of order $p^e$. By above lemma $k(1) = e-2$ and the inequality follows. Conversely, when the inequality holds, we have $||G^{\ast}_1|| = \{ p^e \}$. As $G_1$ is regular, we have ${\mathrm{exp}}(G_1) = p^e = {\mathrm{exp}}(G)$. 

(4) Immediate from (1)-(3). \QED

\begin{definition}
\label{edG-def} Let $G$ be a $p$-group of maximal class of order $p^n \geq p^{p+1}$ and exponent $p^e \geq p^2$ and $\lambda \geq 1$. Then $G$ is said to be of {\bf exponent depth} $ed(G) = \lambda$ if ${\mathrm{exp}}(G_{\lambda}) = p^e$ and ${\mathrm{exp}}(G_{\lambda+1}) \leq p^{e-1}$. 
\end{definition}

\begin{remark}\label{exp-descend} $ed(G) \leq n-2$ as the center $Z(G) = G_{n-1}$ is cyclic of order $p$. Also using {\ref{p+1-lemma}} the conditions $e \geq 2$ and $ed(G) \geq 2$ forces $n \geq p+2$. Hence if $G$ has order $p^n \geq p^{p+2}$ and exponent $p^e \geq p^2$ then 
\[
2 + (e-1)(p-1) \leq \log_p(|G|) \leq 1 + e(p-1)
\]
and $ed(G) = 1$ if and only if $\log_p(|G|) = 2 + (e-1)(p-1)$. Thus for the remaining values of $n$ satisfying the inequality we must have $ed(G) = 2$. This also shows for each $1 \leq N \leq e$, we have some $i(N)$ so that $||G^{\ast}_{i(N)}|| = \{ p^N \}$. Using regularity of $G_1$ we infer that ${\mathrm{exp}}(G_{i(N)}) = p^N$.
\end{remark}


\section{\bf $z$-classes of $p$-groups of maximal class}
\label{z-classes}

The notion of a $z$-class in a group is a weaker compared to the notion of conjugacy class. In general the study of $z$-classes stems from the fact that while the number of conjugacy classes are too many, which in fact is true for finite $p$-groups, the number of $z$-classes are somewhat controlled. Recently in (\cite[Thm. A2]{kkj}) it was proved that the number of $z$-classes in a finite non-abelian $p$-group $G$ is the least if and only if $G$ is of maximal class (up to isoclinism) with $G_1$ abelian, and this least number is $p+2$. Our main interest to develop this section is to get some clarity on the structures of the $z$-classes of the uniform elements beyond the results derived in {\cite{kkj}} that will be useful in later sections. We will start from some definitions. 

\begin{definition} Let $G$ be a group. Two elements $g, h \in G$ are called {\bf $z$-equivalent} if there is an $x \in G$ such that $x^{-1} C_G(g) x = C_G(h)$, i.e., the centralizer subgroups of $g$ and $h$ in $G$ are conjugate in $G$. Whenever this happens, we denote it by $g \sim_z h$ and say that $g$ and $h$ are {\bf $z$-equivalent} in $G$.  
\end{definition}

It is easy to check that $\sim_z$ is an equivalence relation and the its orbits are called the {\bf $z$-classes} of $G$. In fact the $z$-classes are union of some conjugacy classes in $G$. 

\begin{theorem}\label{z-class-unifelts} Let $G$ be a $p$-group of maximal class of order $p^n \geq p^4$. Then:
\begin{enumerate}[(i)]
\item If $c(G) \geq 1$, then the normal set $G^{\ast}_0$ contain $p$ number of $z$-classes.

\item If $c(G) = 0$, then the normal set $G \setminus (G_1 \cup C_G(G_{n-2}))$ contain $p-1$ number of $z$-classes.

\item In either case, for $s \in G \setminus (G_1 \cup C_G(G_{n-2}))$ the $z$-class of $s$ is given by $\cup_{j=1}^{p-1} s^j G_2$. Moreover all the elements in a fixed $z$-class has same order.
\end{enumerate}
\end{theorem}

\bigskip

\noindent {\bf Proof.} Set $T := G \setminus (G_1 \cup C_G(G_{n-2}))$. For any $s \in T$, we know that $s^G = s G_2$ and $C_G(s) = \langle s \rangle Z(G)$ with $|C_G(s)| = p^2$. Two elements $s, s^{\prime} \in T$ are $z$-equivalent if and only if $C_G(s^{\prime}) = g^{-1} C_G(s) g = C_G(g^{-1} sg)$ for some $g \in G$. Since $s^{\prime}, g^{-1} sg \not\in Z(G)$ and $C_G(s^{\prime}) = \langle s^{\prime} \rangle Z(G), C_G(g^{-1} sg) = \langle g^{-1} sg \rangle Z(G)$, and this happens if and only if 
\begin{equation}\label{central-eqn}
g^{-1} sg = (s^{\prime})^{\lambda} z ~{\mathrm{for~some}}~ 1 \leq \lambda \leq p-1, z \in Z(G), g \in G 
\end{equation}
Now fix elements $s \in T, s_1 \in G^{\ast}_1$. If $c(G) \geq 1$, the normal set $T$ has precisely $p(p-1)$ many distinct $G$-conjugacy classes given by
\[
(s^j s_1^{j_1})^G = (s^j s_1^{j_1}) G_2 ~~~~(1 \leq j \leq p-1, 0 \leq j_1 \leq p-1)
\]
Let two such representatives $s^j s_1^{j_1}, s^{j^{\prime}} s_1^{j_1^{\prime}}$ satisfy ({\ref{central-eqn}}) by 
\[
s^j s_1^{j_1} [s^j s_1^{j_1}, g] = g^{-1} (s^j s_1^{j_1}) g = (s^{j^{\prime}} s_1^{j_1^{\prime}})^{\lambda} z
\]
for some $\lambda \in \{ 1, \dotsc, p-1 \}, z \in Z(G), g \in G$. Then $s^j G_1 = s^{j^{\prime} \lambda} G_1$. Hence $j \equiv j^{\prime} \lambda$ mod $p$. Since $s^p \in Z(G) \subseteq G_2$, we have $s_1^{j_1} G_2 = s_1^{j_1^{\prime} \lambda} G_2$. Hence $j_1 \equiv j_1^{\prime} \lambda$ mod $p$. Thus, if two such $G$-conjugacy class combine to form a $z$-conjugacy class, then $(j, j_1) \equiv (j^{\prime}, j_1^{\prime}) \lambda$ mod $p$. 

\smallskip

\noindent Conversely, if for a pair $(j, j_1), (j^{\prime}, j_1^{\prime}) \in {\mathbb F}^{\ast}_p \times {\mathbb F}_p$ which satisfy $(j, j_1) \equiv (j^{\prime}, j_1^{\prime}) \lambda$ mod $p$ for some $\lambda \in {\mathbb F}^{\ast}_p$, we have  
\[
s^j s_1^{j_1} = s^{j^{\prime} \lambda + ap} s_1^{j_1^{\prime} \lambda + bp} = (s^{j^{\prime}} s_1^{j_1^{\prime}})^{\lambda} w
\]
for some integers $a, b$ and some $w \in G_2$, since $s^p \in Z(G)$ and $s_1^p \in G_2$. Now $s^j s_1^{j_1} w^{-1} \in s^j s_1^{j_1} G_2 = (s^j s_1^{j_1})^G$. Hence there exists some $g \in G$ such that 
\[
g^{-1} (s^j s_1^{j_1}) g = s^{j^{\prime}} s_1^{j_1^{\prime}}
\] 
This proves $s^j s_1^{j_1}$ and $s^{j^{\prime}} s_1^{j_1^{\prime}}$ are $z$-conjugate in $G$. Now, since we have $p-1$ many such multiples $\lambda$, these together will form ${\frac {p(p-1)}{p-1}} = p$ many $z$-classes. The case $c(G)=0$ is similar except for the fact that the normal set $T$ contain $(p-1)^2$ many conjugacy classes. This proves everything except the last comment in (iii).

\smallskip 

\noindent From above we have an ${\mathbb F}^{\ast}_p$-action to the uniform elements modulo $G_2$. In particular $g_1, g_2 \in G^{\ast}_0$ belong to the same $z$-class if and only if there exists $1 \leq \lambda \leq p-1$ and $z \in G_2$ so that $g_2 = g_1^{\lambda} z$. Since two elements modulo $G_2$ are conjugate to each other, the last remark follows. \QED

\bigskip

\noindent It is proved in {\cite{kkj}}, that while $G_1$ is abelian, $G$ has precisely $p+2$ conjugacy classes. The following statement now give the precise description of these classes.

\bigskip

\begin{corollary} Let $G$ be a $p$-group of maximal class of order $p^n \geq p^4$ with $G_1$ abelian. Then the $z$-classes of $G$ are given by $Z(G), G_1 \setminus Z(G)$ and the $p$ number of $z$-classes in above theorem.
\end{corollary}

\bigskip

\noindent {\bf Proof.} Since $C_G(a) = G$ if and only if $a \in Z(G)$, the subgroup $Z(G)$ form a single $z$-class. Now let $a \in G_1 \setminus Z(G)$. Then $C_G(a) \supseteq G_1$. Since $a \not\in Z(G)$, this must be equal. This shows $C_G(a) = G_1$ for every $a \in G_1 \setminus Z(G)$. Hence $G_1 \setminus Z(G)$ form another $z$-class.  \QED 

\begin{definition}\label{def-excep} Let $G$ be a $p$-group of maximal class of order $p^n$ with $||G^{\ast}_0|| = \{ p, p^2 \}$ and $\epsilon \in \{ 1, 2 \}$. Then $G$ is said to be of {\bf $p^{\epsilon}$-exceptional type} if there exists an unique $z$-class in $G \setminus (G_1 \cup C_G(G_{n-2}))$ consists of elements of order $p^{\epsilon}$. Otherwise, we say that $G$ is {\bf not of exceptional type}.
\end{definition}

\begin{remark}\label{3^2-excep} In case $|G| \geq p^{p+2}$, the definition of $p$-exceptional (resp. $p^2$-exceptional) simply means if $s \in G^{\ast}_0$ with $|s| = p$ (resp. $|s| = p^2$), then for any $s_1 \in G^{\ast}_1$, we must have $|ss_1| = p^2$ (resp. $|ss_1| = p$).   
\end{remark} 

\noindent Here are some interesting consequences of $p$- or $p^2$-exceptional type towards $(p^a,p^b,p^c)$-generation of $G$.

\begin{definition} A $p$-group $G$ is said to have a $(p^a,p^b,p^c)${\bf{-triple}} $(u,v,w) \in G^3$ if $|u| = p^a, |v| = p^b, |w| = p^c$ and $uvw = 1$. A $(p^a,p^b,p^c)$-triple $(u,v,w)$ of $G$ is called a {\bf {generating triple}} if $G = \langle u, v \rangle$. 
\end{definition}

\begin{remark}\label{p1p2pe-gentriple} Whenever $||G^{\ast}_0|| = \{ p, p^2 \}$ the group $G$ has $(p,p^2,p^e)$ generating triple, i.e., \newline $(0; 1, 1, 0, \dotsc, 0, 1) \in {\mathcal{S}}(G)$ by choosing $s \in G^{\ast}_0, s_1 \in G^{\ast}_1$ with $|s| = p, |ss_1| = p^2, |s_1| = p^e$ and considering the long relation $s^{-1} \cdot (ss_1) \cdot s^{-1}_1 =1$. Also similar choices can be made to show $G$ has both $(p,p,p^e)$ and $(p^2, p^2, p^e)$-generating triple while $G$ is not exceptional. Similarly while $||G^{\ast}_0|| = \{ p^{\epsilon} \}$ with $\epsilon = 1, 2$, $G$ has $(p^{\epsilon},p^{\epsilon},p^e)$-generating triple. 
\end{remark}

\begin{lemma} Let $G$ be a $p$-group of maximal class of order $p^n \geq p^{p+2}$ and ${\mathrm{exp}}(G) = p^e$. Suppose $||G^{\ast}_0|| = \{ p, p^2 \}$.

\begin{enumerate}[(i)] 

\item If $G$ is of $p$-exceptional type, it cannot have $(p,p,p^e)$ triple.

\item If $G$ is of $p^2$-exceptional type, it cannot have $(p^2,p^2,p^e)$ triple. 
\end{enumerate}
\end{lemma}

\noindent {\bf Proof.} (i) Suppose $G$ has such a triple realized as $|u| = p, |v| = p, |\xi| = p^e, uv \xi = 1$. Since $G_1$ is regular, one of $u, v$ must be from $G^{\ast}_0$. Without loss of generality assume $v \in G^{\ast}_0$. Since $\xi \in G^{\ast}_1$ and $G$ is of $p$-exceptional type, $|v \xi| = p^2$, which is a contradiction. Proof of (ii) is similar. \QED   

 
\section{\bf Results on groups of exponent depth $\geq 2$}\label{resexp3ed2}

\noindent For this section we will always assume $G$ is a $p$-group of maximal class of order $p^n \geq p^{p+2}$, with exp$(G) = p^e \geq p^2$ and $ed(G) = \lambda \geq 2$. Using {\ref{exp-p-p2-lemma}(3) we have $||G^{\ast}_1|| = \{ p^e \}$. We first derive the results while $e \geq 3$ (by {\ref{exp-descend}}, this mean $n \geq 2p$). The situation with $e = 2$ and $ed(G) \geq 2$ is similar (see {\ref{signature-ed2-exp-p2}}). We first start with the case $h \geq 2$.

\begin{theorem}\label{orbgen1or2} Let $G$ be a $p$-group of maximal class of order $p^n \geq p^{p+2}$, ${\mathrm{exp}}(G) = p^e \geq p^3$ and $ed(G) \geq 2$. 

\begin{enumerate}[(a)]

\item If $h \geq 2$ and $(m_1, \dotsc, m_e) \in {\mathbb N}^{e}_{\geq 0}$, then $(h; m_1, \dotsc, m_e) \in {\mathcal S}(G)$.
\item If $h=1$ and $(m_1, \dotsc, m_e) \in {\mathbb N}^{e}_{\geq 0}$ with $m_e \geq 1$, then $(1; m_1, \dotsc, m_e) \in {\mathcal S}(G)$.
\end{enumerate}
\end{theorem}

\noindent {\bf Proof.} (a) Let $s \in G^{\ast}_0, s_1 \in G^{\ast}_1$. Then $G = \langle s, s_1 \rangle$. Let $N$ be the largest integer among $\{ 1, 2, \dotsc, e \}$, such that $m_N \neq 0$. Since $ed(G) \geq 2$ we have $||G^{\ast}_1|| = \{ p^e \}$ and consequently by {\ref{exp-descend}} and {\ref{exp-p-p2-lemma}}(3), $||G^{\ast}_{n - N(p-1)}|| = \{ p^N \}$. Set $j := n - N(p-1)$. The element $s_j = [s_{j-1},s] \in G_j$ has order $p^N$. Now set $y = [s_{j-1}, s]$ and we have a relation $[s_1, s_1] [s, s_{j-1}] y = 1$. Now using A2.2, we have (a).

(b) For any $s \in G^{\ast}_0, s_1 \in G^{\ast}_1$ with $|[s, s_1]| = p^e$. Set $y := [s_1, s]$. Then we have the relation $[s, s_1] y = 1$ in $G$ where $y$ has order $p^e$. 
Now (b) follows by A2.2. \QED

\smallskip

Note that every uniform element $s \in G^{\ast}_0$ satisfy $s^p \in Z(G)$. Hence $||G^{\ast}_0|| \subseteq \{ p, p^2 \}$. Thus we have three possible cases : $||G^{\ast}_0||$ can be $\{ p \}, \{ p^2 \}$, or $\{ p, p^2 \}$. The following result is useful for both cases ${\mathrm{ed}}(G) = 1$ and $\geq 2$.

\begin{proposition}\label{gen1} Let $G$ be a $p$-group of maximal class of order $p^n \geq p^{p+2}$, ${\mathrm{exp}}(G) = p^e \geq p^3$. Suppose $\sigma = (1;m_1, \dotsc,m_e) \in {\mathcal{S}}(G)$ realized by a generating set $\Lambda$ as in A2.1. 

\noindent If $1 \leq N(\sigma) \leq e-1$, then $|G^{\ast}_0 \cap (E_1(\Lambda) \cup E_2(\Lambda))| \geq 2$. In case of equality either $|G^{\ast}_0 \cap E_1(\Lambda)| = 2$ or $|G^{\ast}_0 \cap E_2(\Lambda)| = 2$.

\end{proposition}   

\noindent {\bf Proof.} (a) Using Burnside basis theorem either of the following occur : (i) $G = \langle a_1, b_1 \rangle$, (ii) $G = \langle a_1, X_{i_1 j_1} \rangle$ or $G = \langle b_1, X_{i_1 j_1} \rangle$ with $i_1 \in \{ 1, 2 \}$, (iii) $a_1, b_1 \in G_2$ and $G = \langle X_{i_1 j_1}, X_{i_2 j_2} \rangle$ with $i_1, i_2 \in \{ 1, 2 \}$.

\noindent In case (i), from {\ref{exp-p-p2-lemma}}(3) we have $|[a_1, b_1]| = p^e$ and using commutator identities we may assume $(a_1, b_1) \in G^{\ast}_0 \times G^{\ast}_1$. If all the elliptic generators $X_{ij} \in G_1$, this contradict $N < e$ using regularity of $G_1$. Hence some of $X_{ij}$ must be from $G^{\ast}_0$. As $G_1 \unlhd G$ maximal subgroup, using the long relation we have at least two of $X_{ij}$ belong $G^{\ast}_0$. Now assume $G^{\ast}_0 \cap (E_1(\Lambda) \cup E_2(\Lambda)) = \{ X_{i_1 j_1}, X_{i_2 j_2} \}$ where $i_1, i_2 \in \{ 1, 2 \}$. We claim $i_1 = i_2$ : if this is not true, we may assume $i_1 = 1, i_2 = 2$. Combining the remaining generators of long relation we reduce to the equation $X_{i_1 j_1} X_{i_2 j_2} y = 1$ for some $y \in G_2$. Using {\ref{unif-prop}}, this means $X^{-1}_{i_1 j_1}$ is a conjugate of $X_{i_2 j_2}$, a contradiction.

\noindent In case (ii), assume $\langle a_1, b_1 \rangle \neq G = \langle a_1, X_{i_1 j_1} \rangle$. Since $\{ a_1, b_1, X_{i_1 j_1} \}$ is ${\mathbb{F}}_p$-linearly dependent modulo $G_2$, using commutator identities we may assume $(a_1, b_1) \in G^{\ast}_0 \times G_2$. Since $||G^{\ast}_1|| = \{ p^e \}$ we have $X_{i_1 j_1} \in G^{\ast}_0$. As above $|G^{\ast}_0 \cap (E_1(\Lambda) \cup E_2(\Lambda))| \geq 2$ and the proof is similar. 

\noindent In situation (iii), if $|G^{\ast}_0 \cap (E_1(\Lambda) \cup E_2(\Lambda))| = 2$, then as above $X_{i_1 j_1}$ is conjugate to $X_{i_2 j_2}$. But then these are ${\mathbb{F}}_p$-linearly dependent modulo $G_2 = \Phi(G)$, a contradiction.  \QED 

\begin{theorem}\label{h1me0ed2} Let $G$ be a $p$-group of maximal class of order $p^n \geq p^{p+2}$, ${\mathrm{exp}}(G) = p^e \geq p^3$ and $ed(G) \geq 2$. Let $2 \leq N \leq e-1$ and consider $\sigma = (1; m_1, \dotsc, m_N, 0, \dotsc, 0) \in {\mathbb N}^{e+1}_{\geq 0}$ with $m_N \geq 1$.

\begin{enumerate}[(i)]

\item While $||G^{\ast}_0|| = \{ p \}$, we have $\sigma \in {\mathcal S}(G)$ if and only if $m_1 \geq 2$. 

\item While $||G^{\ast}_0|| = \{ p^2 \}$, we have $\sigma \in {\mathcal S}(G)$ if and only if $m_2 \geq 2$.

\item While $||G^{\ast}_0|| = \{ p, p^2 \}$, we have $\sigma \in {\mathcal S}(G)$ if and only if either $m_1 \geq 2$ or, $m_2 \geq 2$.

In case $N = 1$ we have the following :

\item while $||G^{\ast}_0|| = \{ p \}$, or $\{ p, p^2 \}$ we have $(1; m_1, 0, \dotsc, 0) \in {\mathcal S}(G)$ if and only if $m_1 \geq 2$. 

\item while $||G^{\ast}_0|| = \{ p^2 \}$, we have $(1; m_1, 0, \dotsc, 0) \not\in {\mathcal S}(G)$. 
\end{enumerate}
\end{theorem}
 
\noindent {\bf Proof.} The necessity of the statements (i)-(iv) follows from {\ref{gen1}}. If $N=1$ and $||G^{\ast}_0|| = \{ p^2 \}$, then from {\ref{gen1}} $m_2 \neq 0$ is necessary as well. We are left to show that for (i)-(iv), these are sufficient. 

\smallskip

\noindent {\it Case-I: $N \geq 3$}: In each case pick a pair $s \in G^{\ast}_0, s_1 \in G^{\ast}_1$ and an element $x \in G_2$ with $|x| = p^N$. Then the set $\Lambda = \{ a_1 = s_1, b_1 = s_1, x_1 = s^{-1}, x_2 = sx^{-1}, y = x \}$ generate $G$ and satisfy the long relation. In case (i) $|s| = p$ and $x_2 = sx^{-1}$ is conjugate to $s$ by {\ref{unif-prop}} and hence $|x_2| = p$. Now using \ref{basic} we have $(1; m_1, \dotsc, m_N, 0, \dotsc, 0) \in {\mathcal{S}}(G)$. The case (ii) is similar with $|x_1| = |x_2| = p^2$. In case (iii) we may choose $s$ with either $|s| = p$ or $|s| = p^2$. 

\smallskip

\noindent {\it Case-II: $N = 1, 2$}: If $||G^{\ast}_0|| = \{ p^2 \}$, the generating set $\Lambda^{\prime} = \{ a_1 = s_1, b_1 = s_1, x_1 = s, x_2 = s^{-1} \}$ can be extended to obtain the corresponding $\sigma \in {\mathcal {S}}(G)$ with the minimal condition $m_2 \geq 2$. If $||G^{\ast}_0|| = \{ p \}$, the above set $\Lambda$ with minimal condition $m_1 \geq 2$ works. While $||G^{\ast}_0|| = \{ p, p^2 \}$, we require either $m_2 \geq 2$ or $m_1 \geq 2, m_2 \geq 1$ (as $m_2 \neq 0$) and obtain the corresponding generating sets $\Lambda^{\prime}$ or $\Lambda$ respectively. Finally when $N=1$, and $p \in ||G^{\ast}_0||$ we simply use $\Lambda^{\prime}$ as above. \QED

\smallskip

The following results do not use the condition ${\mathrm{ed}}(G) \geq 2$, and hence could be used if we ensure ${\mathrm{exp}}(G) \geq p^3$. This will be useful for the case ${\mathrm{ed}}(G) = 1$ in the next section.  

\begin{theorem}\label{0-genus-m_e-one} Let $G$ be a $p$-group of maximal class of order $p^n \geq p^{p+2}$ and ${\mathrm{exp}}(G) = p^e \geq p^3$. Consider the tuple $\sigma = (0; m_1, m_2, \dotsc, m_e)$ with $m_e \geq 1$. 
\begin{enumerate}[(i)]
\item If $||G^{\ast}_0|| = \{ p \}$, we have $\sigma \in {\mathcal S}(G)$ if and only if $m_1 \geq 2$.

\item If $||G^{\ast}_0|| = \{ p^2 \}$, we have $\sigma \in {\mathcal S}(G)$ if and only if $m_2 \geq 2$.

\item If $||G^{\ast}_0|| = \{ p, p^2 \}$, we have $\sigma \in {\mathcal S}(G)$ if and only if either $m_1 \geq 1, m_2 \geq 1$ or one of the following conditions hold:
\end{enumerate}
\begin{enumerate}[(a)]
\item $m_2 \geq 2$ if $G$ is of $p$-exceptional type,

\item $m_1 \geq 2$ if $G$ is of $p^2$-exceptional type,

\item either $m_1 \geq 2$ or $m_2 \geq 2$ if $G$ is not of exceptional type.
\end{enumerate} 
\end{theorem} 

\noindent {\bf Proof.} The sufficiency follows from the definition {\ref{def-excep}}. We need to prove the necessity of the conditions. Let $\sigma \in {\mathcal{S}}(G)$ is realized by a generating set $\Lambda$ as in A2.1. As seen before we must have at least two elliptic generators $u, v \in G^{\ast}_0$ with $G = \langle u, v \rangle$. In the first two cases we must have $|u| = p^{\epsilon} = |v|$ where $||G^{\ast}_0|| = \{ p^{\epsilon} \}$. In case (iii) while $|u| = p, |v| = p^2$ we see that $m_1 \geq 1, m_2 \geq 1$. Also, while $G$ is not of exceptional type, using {\ref{def-excep}}, we have either $m_1 \geq 2$ or $m_2 \geq 2$. In case $G$ is of $p$-exceptional type, two generator of order $p$ cannot generate $G$, by {\ref{def-excep}}. Thus either $m_1 \geq 1, m_2 \geq 1$, or else $m_2 \geq 2$. Similarly while $G$ is of $p^2$-exceptional type, we have either $m_1 \geq 1, m_2 \geq 1$, or else $m_1 \geq 2$. This completes the proof.
\QED

\begin{table}[h]
  \begin{center}
  \caption {${\mathrm{ed}}(G) \geq 2$, ${\mathrm{exp}}(G) \geq p^3$} \label{ed2table}
\begin{tabular}{c   c   c   c}
\hline
$h$ & $N$ & $m_N$ & property of $G$ \\
\hline
$\geq 2$ & $0$      & $0$      &  -  \\
         & $\neq 0$ & - &  -  \\
    $1$  & $e$      & - &  -  \\
         & $2 \leq N \leq e-1$    & $m_1 \geq 2$ & $||G^{\ast}_0|| = \{ p \}$ \\ 
         &          & $m_2 \geq 2$ & $||G^{\ast}_0|| = \{ p^2 \}$ \\ 
         &          & $m_1 \geq 2$ or $m_2 \geq 2$ & $||G^{\ast}_0|| = \{ p, p^2 \}$ \\ 
         & $1$    & $m_1 \geq 2$ & $p \in ||G^{\ast}_0||$ \\
    $0$  & $e$    & $m_1 \geq 2$ & $||G^{\ast}_0|| = \{ p \}$  \\         
    		 & 		 & $m_2 \geq 2$ & $||G^{\ast}_0|| = \{ p^2 \}$  \\         
    		 & 		 & $m_1 + m_2 \geq 2$ & $||G^{\ast}_0|| = \{ p, p^2 \}$, not exceptional \\ 		 & 		 & $m_1 \geq 1, m_2 \geq 1$ or $m_2 \geq 2$ & $||G^{\ast}_0|| = \{ p, p^2 \}$, $p$-exceptional \\
    		 & 		 & $m_1 \geq 1, m_2 \geq 1$ or $m_1 \geq 2$ & $||G^{\ast}_0|| = \{ p, p^2 \}$, $p^2$-exceptional \\ 
    		 & $N \leq e-1$ & $m_1 \geq 3$ & $||G^{\ast}_0|| = \{ p \}$  \\   
    		 & 		 & $m_2 \geq 3$ & $||G^{\ast}_0|| = \{ p^2 \}$  \\
    		 & 		 & $m_1 + m_2 \geq 3$ & $||G^{\ast}_0|| = \{ p, p^2 \}$, not exceptional  \\
    		 & 		 & $m_1 \geq 1, m_2 \geq 2$ & $||G^{\ast}_0|| = \{ p, p^2 \}$, $p$-exceptional  \\
    		 & 		 & $m_1 \geq 2, m_2 \geq 1$ & $||G^{\ast}_0|| = \{ p, p^2 \}$, $p^2$-exceptional  \\       
\hline
\end{tabular}
\end{center}
\end{table}

\begin{theorem}\label{0-genus-m_e-0} Let $G$ be a $p$-group of maximal class of order $p^n$ and ${\mathrm{exp}}(G) = p^e \geq p^3$. Consider $\sigma = (0; m_1, m_2, \dotsc, m_N, 0, \dotsc, 0) \in {\mathbb N}^{e+1}_{\geq 0}$ with $m_N \geq 1$ and $N \leq e-1$. 
\begin{enumerate}[(i)]
\item If $||G^{\ast}_0|| = \{ p \}$, we have $\sigma \in {\mathcal S}(G)$ if and only if $m_1 \geq 3$.

\item If $||G^{\ast}_0|| = \{ p^2 \}$, we have $\sigma \in {\mathcal S}(G)$ if and only if $m_2 \geq 3$.

\item If $||G^{\ast}_0|| = \{ p, p^2 \}$, we have $\sigma \in {\mathcal S}(G)$ if and only if $(m_1, m_2)$ satisfy:
\end{enumerate}
\begin{enumerate}[(a)]
\item $m_1 \geq 1, m_2 \geq 2$ if $G$ is of $p$-exceptional type,

\item $m_1 \geq 2, m_2 \geq 1$ if $G$ is of $p^2$-exceptional type,

\item either of the conditions in (a) or (b) if $G$ is not of exceptional type. 
\end{enumerate}
\end{theorem}

\noindent {\bf Proof.} (Necessity) As seen in the proof of {\ref{0-genus-m_e-one}} we require two elliptic generators $u, v \in G^{\ast}_0$ with $G = \langle u, v \rangle$. Since $N \leq e-1$, the elliptic generators belong to $G^{\ast}_0 \cup G_2$. Since $u, v$ are ${\mathbb{F}}_p$-linearly independent modulo $G_2$, from long relation we have at least three of them from $G^{\ast}_0$. The remaining argument of the necessary part is similar as in the previous theorem. 

\noindent For sufficiency while $N > 2$, let $x \in G_2$ with $|x| = p^N$. Now choose $s \in G^{\ast}_0, s_1 \in G^{\ast}_1$ and set $\Lambda = \{ x_1 = s, x_2 = ss_1, x_3 = (s^2 s_1)^{-1} x^{-1}, x_4 = x \}$, which satisfy the criterion. In case $N = 1, 2$ the generating set $\Lambda = \{ x_1 = s, x_2 = ss_1, x_3 = (s^2 s_1)^{-1} \}$ realize $\sigma \in {\mathcal S}(G)$. \QED 

Now we discuss the case ${\mathrm{exp}}(G) = p^2$ and ${\mathrm{ed}}(G) \geq 2$. From \ref{exp-descend} it follows that $p^{p+2} \leq |G| \leq p^{2p-1}$. In this case the groups which are not of exceptional type require some special attention. 

\begin{definition} A finite $p$-group is said to be of $(p,p,p)$-type (resp. $(p^2,p^2,p^2)$-type) if there are elements $x, y \in G \setminus \{ 1 \}$ so that $|x| = |y| = |xy| = p$ (resp. $|x| = |y| = |xy| = p^2$) and $G = \langle x, y \rangle$. 
\end{definition}

\noindent A finite $p$-group of maximal class of order $p^n$ is of $(p,p,p)$-type if and only if $G = \langle s, {\tilde{s}} \rangle$ so that the distinct maximal subgroups $M_{s} := \langle s, G_2 \rangle, M_{{\tilde{s}}} = \langle {\tilde{s}}, G_2 \rangle$ are different from $G_1, C_G(G_{n-2})$ and $||M_{s} \setminus G_2|| = \{ p \} = ||M_{{\tilde{s}}} \setminus G_2||$. In other words, $s$ and ${\tilde{s}}$ belong to distinct $z$-classes of $G$ with elements of order $p$. The following result connects the idea of $(p,p,p)$-type in our setting.

\begin{lemma}\label{(p,p,p)-type-lem} Let $G$ be a $p$-group of maximal class of order $p^n$ with $p+2 \leq n \leq 2p-1$. Let $\kappa(G)$ denote the number of $z$-classes in $G^{\ast}_0$ which contain elements of order $p$.
\begin{enumerate}[(i)]
\item If $||G^{\ast}_0|| = \{ p \}$, then $G$ is of $(p,p,p)$-type.
\item If $G$ is $p^2$-exceptional with $p \geq 5$, then $G$ is of $(p,p,p)$-type.
\item If $0 \leq \kappa(G) \leq 2$ (this includes $||G^{\ast}_0|| = \{ p^2 \}$, $G$ is $p$-exceptional, and all $3$-groups except $||G^{\ast}_0|| = \{ 3 \}$), then $G$ is not of $(p,p,p)$-type.
\item If $||G^{\ast}_0|| = \{ p, p^2 \}$, $G$ is not exceptional with $\kappa(G) \geq 3$, then $G$ can be of both types $(p,p,p)$ or otherwise. 
\end{enumerate}
\end{lemma}  

\smallskip

\noindent {\bf Proof :} By hypothesis and \ref{exp-descend}, we have ${\mathrm{ed}}(G) \geq 2$. If $G$ is of $(p,p,p)$-type with $G = \langle x, y \rangle$ and $x^p = y^p = (xy)^p = 1$, then all three elements $x, y, xy$ must belong to distinct $z$-classes consisting of elements of order $p$ as $||G^{\ast}_1 \cup G^{\ast}_2|| = \{ p^2 \}$. By \ref{z-class-unifelts}(iii), the statement (iii) follows. For a choice of $(s, s_1) \in G^{\ast}_0 \times G^{\ast}_1$, setting $x = s, y = ss_1$, the statement (i) follows. 

\smallskip

\noindent Proof of (ii) : start with a choice of $(s, s_1) \in G^{\ast}_0 \times G^{\ast}_1$. Adopting notations in the proof of \ref{z-class-unifelts}, the $p$ many $z$-classes (since $c(G) \geq 1$) in $G^{\ast}_0$ are represented by elements $\{ s \equiv (1,0), ss_1 \equiv (1,1), \dotsc, ss^{p-1}_1 \equiv (1,p-1) \}$ (the co-ordinate tuples represent elements in ${\mathbb F}^{\ast}_p \times {\mathbb F}_p$ modulo ${\mathbb F}^{\ast}_p$ as in \ref{z-class-unifelts}). Now assuming $|s| = p^2$, set $x = ss_1, y = ss^2_1$. Then $xy \equiv s^2 s^3_1 ~{\mathrm{mod}} G_2 \equiv (2,3) \not\equiv (1,0)$ modulo ${\mathbb F}^{\ast}_p$ (since $p \geq 5$). Hence $xy \not\in$ the $z$-class containing $s$, i.e., $|xy| = p$. In (iv) $G$ has $(p,p,p)$-type if and only if there are two $z$-classes with order $p$ elements represented by $(1, i_1), (1, i_2)$ so that $(2, i_1 + i_2) \equiv (1, \lambda(i_1 + i_2))$ represents another $z$-class with elements of order $p$, where $\lambda$ is the inverse of $2$ in ${\mathbb F}^{\ast}_p$. \QED  

\begin{table}[h]
  \begin{center}
  \caption {${\mathrm{ed}}(G) \geq 2$, ${\mathrm{exp}}(G) = p^2$} \label{ed2table-expp2}
\begin{tabular}{c   c   c   c}
\hline
$h$ & $N$ & $m_N$ & property of $G$ \\
\hline
$\geq 2$ & $0$      & $0$      &  -  \\
         & $\neq 0$ & - &  -  \\
    $1$  & $2$      & - &  -  \\
         & $1$    & $m_1 \geq 2$ & $p \in ||G^{\ast}_0||$ \\ 
    $0$  & $2$    & $m_1 \geq 2$ & $||G^{\ast}_0|| = \{ p \}$  \\         
    		 & 		 & $m_2 \geq 3$ & $||G^{\ast}_0|| = \{ p^2 \}$  \\         
    		 & 		 & $m_1 + m_2 \geq 3$ & $||G^{\ast}_0|| = \{ p, p^2 \}$, not exceptional \\
    		 & 		 & $m_2 \geq 2, m_1 + m_2 \geq 3$ & $||G^{\ast}_0|| = \{ p, p^2 \}$, $p$-exceptional \\
    		 & 		 & $m_1 \geq 2, m_1 + m_2 \geq 3$ & $||G^{\ast}_0|| = \{ p, p^2 \}$, $p^2$-exceptional \\ 
    		 & $1$ & $m_1 \geq 3$ & $||G^{\ast}_0|| = \{ p \}$, or $G$ is $p^2$-exceptional with $p \geq 5$  \\ 
    		 &       &    &  $||G^{\ast}_0|| = \{ p, p^2 \}$, $p \geq 5$, $G$ not exceptional, $(p,p,p)$-type \\  
    		 & 		 & $m_1 \geq 4$ & $||G^{\ast}_0|| = \{ p, p^2 \}$, $p \geq 5$, $G$ not exceptional, not $(p,p,p)$-type \\
    		 &		 & 				& $||G^{\ast}_0|| = \{ 3, 3^2 \}$, $G$ is $3^2$-exceptional. \\
\hline    		 
\end{tabular}
\end{center}
\end{table}

\begin{theorem}\label{signature-ed2-exp-p2} Let $G$ be a finite $p$-group of maximal class of order $p^n \geq p^{p+2}$ and exponent ${\mathrm{exp}}(G) = p^2$. Suppose that $G$ has exponent depth ${\mathrm{ed}}(G) \geq 2$. Then the following statements hold for $\sigma = (h; m_1, m_2)$ :

\begin{enumerate}[(i)]

\item If $h \geq 2$, then for any $(m_1, m_2) \in {\mathbb N}^{2}_{\geq 0}$, we have $(h; m_1, m_2) \in {\mathcal S}(G)$.
\item If $h=1$, then for any $(m_1, m_2) \in {\mathbb N}^{2}_{\geq 0}$ with $m_2 \geq 1$, we have $(1; m_1, m_2) \in {\mathcal S}(G)$.
\item $(1; m_1, 0) \in {\mathcal{S}}(G)$ if and only if $m_1 \geq 2$ and $p \in ||G^{\ast}_0||$.

\item For $\sigma = (0; m_1, m_2)$ with $m_2 \geq 1$ and $m_1 + m_2 \geq 3$, we have 

\begin{enumerate}[(a)]

\item If $||G^{\ast}_0|| = \{ p \}$, $\sigma \in {\mathcal{S}}(G)$ if and only if $m_1 \geq 2$.

\item If $||G^{\ast}_0|| = \{ p^2 \}$, $\sigma \in {\mathcal{S}}(G)$ if and only if $m_2 \geq 3$.

\item If $||G^{\ast}_0|| = \{ p, p^2 \}$ and $G$ not exceptional, $\sigma \in {\mathcal{S}}(G)$ always.

\item If $||G^{\ast}_0|| = \{ p, p^2 \}$ and $G$ is $p$-exceptional, $\sigma \in {\mathcal{S}}(G)$ if and only if $m_1 + m_2 \geq 3$ and $m_2 \geq 2$.

\item If $||G^{\ast}_0|| = \{ p, p^2 \}$ and $G$ is $p^2$-exceptional, $\sigma \in {\mathcal{S}}(G)$ if and only if $m_1 + m_2 \geq 3$ and $m_1 \geq 2$.

\end{enumerate}

\item For $\sigma = (0; m_1, 0)$ with $m_1 \geq 3$, we have 

\begin{enumerate}

\item If $||G^{\ast}_0|| = \{ p \}$, then $\sigma \in {\mathcal{S}}(G)$ if and only if $m_1 \geq 3$.

\item If either $||G^{\ast}_0|| = \{ p^2 \}$ or $G$ is $p$-exceptional, then $\sigma \not\in {\mathcal S}(G)$.

\item If $||G^{\ast}_0|| = \{ p, p^2 \}$, $G$ is $p^2$-exceptional, then $\sigma \in {\mathcal{S}}(G)$ if and only if $p \geq 5$ and $m_1 \geq 3$ (resp. $p = 3$ and $m_1 \geq 4$).   

\item If $p \geq 5$, $||G^{\ast}_0|| = \{ p, p^2 \}$, $G$ is not exceptional and is of $(p,p,p)$-type (resp. not of $(p,p,p)$-type), then $\sigma \in {\mathcal{S}}(G)$ if and only if $m_1 \geq 3$ (resp. $m_1 \geq 4$).   

\end{enumerate}

\end{enumerate}
\end{theorem}

\noindent {\bf Proof :} Proof of (i) and (ii) is same as \ref{orbgen1or2}. For (iii), if $m_1 = 1$, the generating set $\Lambda = \{ a_1, b_1, X_{11} \}$ satisfy $G = \langle a_1, b_1 \rangle$, using the long relation. This implies $[a_1, b_1] \in G^{\ast}_2$ and hence $|[a_1, b_1]| = p^2$, a contradiction. Hence $m_1 \geq 2$. Now if $||G^{\ast}_0|| = \{ p^2 \}$, we have $||G^{\ast}_0 \cup G^{\ast}_1|| = \{ p^2 \}$ as ${\mathrm{ed}}(G) \geq 2$. Thus each elliptic generator is from $G_2 = \Phi(G)$ forcing $G = \langle a_1, b_1 \rangle$ again. Using regularity of $G_2$, it is not possible. On the other hand while $p \in ||G^{\ast}_0||$ choose $s, s^{\prime} \in G^{\ast}_0, s_1 \in G^{\ast}_1$ with $|s^{\prime}| = p$. Set $\Lambda = \{ a_1 = s, b_1 = s_1, X_{11} = s^{\prime}, X_{12} = ([s, s_1] s^{\prime})^{-1} \}$ and noticing that $X^{-1}_{11}$ and $X_{12}$ are conjugate to each other, the proof follows. 

\noindent Proof of (iv): here we have $h = 0$ and $N = 2$. Hence $m_2 \geq 1$ and $m_1 + m_2 \geq 3$ is necessary for all groups. For (iv)(a) $||G^{\ast}_0|| = \{ p \}$, if $m_1 = 1$, only one elliptic generator is allowed to belong to $G^{\ast}_0$, contradicting the long relation. For $s \in G^{\ast}_0, s_1 \in G^{\ast}_1$, the generating set $\Lambda_1 = \{ X_{11} = s, X_{12} = s^{-1}s_1, X_{21} = s^{-1}_1 \}$ achieve the signature $(0;2,1)$. For (iv)(b) if $m_2 = 2$, at the most two elliptic generators $\in G^{\ast}_0 \cup G^{\ast}_1$ and by long relation they are ${\mathbb F}_p$-linearly dependent, a contradiction. Hence, $m_2 \geq 3$. Now the generating set $\Lambda_1$ with choices $s, s^{\prime}, s_1$ in this case achieve the signature $(0;0,3)$. For (iv)(c) $||G^{\ast}_0|| = \{ p, p^2 \}$ and $G$ is not exceptional (this require $p \geq 5$) it is enough to prove all the three signatures $(0;2,1), (0;1,2), (0;0,3) \in {\mathcal S}(G)$. For this choose $s \in G^{\ast}_0, s_1 \in G^{\ast}_1$ in each of the cases $(|s|, |ss_1|) = (p,p), (p,p^2), (p^2, p^2)$ respectively. Now the generating set $\Lambda_2 = \{ s, s_1, (ss_1)^{-1} \}$ achieve these signatures. Proof of (iv)(d)-(e) is same.  

\noindent Proof of (v): Note that a $3$-group $G$ with $||G^{\ast}_0|| = \{ 3, 3^2 \}$ must be either $3$-exceptional or $3^2$-exceptional by \ref{z-class-unifelts}. Thus we require $p \geq 5$. The proof follows from \ref{(p,p,p)-type-lem} except the second case where we choose $(s, s_1) \in G^{\ast}_0 \times G^{\ast}_1$ with $|s| = p = |ss_1|$ and set $\Lambda_3 = \{ X_{11} = s, X_{12} = ss_1, X_{13} = (ss_1)^{-1}, X_{14} = s^{-1} \}$ which achieve the signature $(0;4,0)$. \QED  

\smallskip

\noindent All the results of this section are summarised in tables (\ref{ed2table}) and (\ref{ed2table-expp2}). 

\section{\bf Results on groups of exponent depth $1$ and $|G| \geq p^{p+2}$}\label{resed1}

\noindent Among the groups of exponent depth one, there are two cases: 
\begin{enumerate}[(i)]
\item groups of order $p^{p+1}$ which has ${\mathrm{exp}}(G) = p^2$, where we may have $c(G) = 0$, which is difficult to analyse compared to the case $c(G) \geq 1$, and
\item groups of order $p^n \geq p^{p+2}$ which always comes with $c(G) \geq 1$. This implies ${\mathrm{exp}}(G) = p^e \geq p^3$ and $n = (e-1)(p-1) + 2$ by \ref{exp-p-p2-lemma}.
\end{enumerate}

\noindent We will discuss the first case in next section. In either situation we have ${\mathrm{exp}}(G_2) = p^{e-1}$ by {\ref{exp-descend}}. 




\begin{remark} In these groups we do not have any particular advantage while we consider $(h; m_1, \dotsc, m_e)$ with $h \geq 2$ and $m_i \neq 0$ for some $i$ against the situation $h \leq 1$ (compare {\ref{orbgen1or2}}). Thus we divide the results based on whether $m_e = 0$ and $m_e \neq 0$. The first result deals with the case $h, m_e \geq 1$.
\end{remark}

\begin{theorem}\label{h2cd1} Let $G$ be a $p$-group of maximal class of order $p^n \geq p^{p+2}$, ${\mathrm{exp}}(G) = p^e$ and ${\mathrm{ed}}(G)=1$. Consider the tuple $\sigma = (h; m_1, \dotsc, m_e)$ with $h \geq 1$ and $m_e \geq 1$. Then:

\begin{enumerate}[(i)]

\item If $m_e \geq 2$, then $\sigma \in {\mathcal{S}}(G)$.

\item If $m_e = 1$, then $\sigma \in {\mathcal{S}}(G)$ if and only if 
\begin{enumerate}[(a)]
\item $m_1 \geq 2$ while $||G^{\ast}_0|| = \{ p \}$,

\item $m_2 \geq 2$ while $||G^{\ast}_0|| = \{ p^2 \}$,

\item $m_1 + m_2 \geq 2$ while $||G^{\ast}_0|| = \{ p, p^2 \}$, $G$ is not exceptional,

\item $m_1 + m_2 \geq 2$ and $m_2 \geq 1$ while $||G^{\ast}_0|| = \{ p, p^2 \}$, $G$ is $p$-exceptional,

\item $m_1 + m_2 \geq 2$ and $m_1 \geq 1$ while $||G^{\ast}_0|| = \{ p, p^2 \}$, $G$ is $p^2$-exceptional.
\end{enumerate}
\end{enumerate} 
\end{theorem}

\noindent {\bf Proof.} (i) This can done using the generating set $\Lambda = \{ a_1 = s = b_1, X_{e1} = s_1, X_{e2} = s^{-1}_1 \}$  
by choosing any $s \in G^{\ast}_0, s_1 \in G^{\ast}_1$ with $|s_1| = p^e$ and by extension principle (A2.2). 

\noindent Necessity of (ii) : Let $\Lambda$ be a generating set of $G$ corresponding to $\sigma = (h; m_1, \dotsc, m_{e-1}, 1)$ with $h \geq 1$. Since $G_1$ is regular and $m_e = 1$ we arrive at a contradiction if the elliptic generator of order $p^e$ is from $G_1$. Thus at least two elliptic generators are from $G^{\ast}_0$. Since $||G^{\ast}_0|| \subseteq \{ p, p^2 \}$ we need $m_1 + m_2 \geq 2$. This immediately proves the necessity of (ii)(a)-(c). 

To prove the necessity of (ii)(d)-(e) we need to check while $G$ is $p$-exceptional, $m_2 = 0$ is not possible (similarly while $G$ is $p^2$-exceptional $m_1 \neq 0$). If $m_2 = 0, m_e = 1$, the long relation reduces to $u_1 \dotsc u_k v = 1$ where $u_i \in G^{\ast}_0, |u_i| = p$ for each $i$ and $|v| = p^e$, since all elements of order $\geq p^3$ are from $G_1$ which is regular. Since $G$ is $p$-exceptional, using {\ref{z-class-unifelts}} there exists a $s \in G^{\ast}_0$ with $|s| = p$ so that each $u_i$ are of the form $u_i = s^{t_i} w_i$ with $1 \leq t_i \leq p-1$ and $w_i \in G_2$. Using commutator identities the above equation reduces to $s^{t_1 + \dotsc + t_k} w^{\prime} = v^{-1}$ for some $w^{\prime} \in G_2$. Since $v \in G^{\ast}_1$, we must have $t_1 + \dotsc + t_k \equiv 0$ mod $p$. But $s^p \in Z(G) \subseteq G_2$. This implies $v \in G_2$, a contradiction. Similar argument works for $p^2$-exceptional case.

\noindent (Sufficiency) For (ii)(a)-(b) choose any $s \in G^{\ast}_0, s_1 \in G^{\ast}_1$ and consider $\Lambda_1 = \{ X_{i1} = s, X_{i2} = s^{-1}s_1, X_{e1} = s^{-1}_1 \}$ where $i=1$ in (a) and $i=2$ in (b). Now these generators can be extended. While $||G^{\ast}_0|| = \{ p, p^2 \}$, the sufficiency of $m_1 + m_2 \geq 1$ follows from {\ref{p1p2pe-gentriple}}. While $G$ is exceptional we can use $\Lambda_1$ with suitable generators and extend if necessary.
\QED

\begin{theorem} Let $G$ be a $p$-group of maximal class of order $p^n \geq p^{p+2}$, ${\mathrm{exp}}(G) = p^e$ and $ed(G)=1$. Consider the tuple $\sigma = (h; m_1, \dotsc, m_N, 0, \dotsc, 0) \in {\mathbb N}^{e+1}_{\geq 0}$ where $0 \leq N = N(\sigma) \leq e-1$ and $h \geq 1$.
\begin{enumerate}[(i)]
\item For $N = 0$, we have $\sigma \in {\mathcal{S}}(G)$ if and only if $h \geq 2$.

\item For $N = 1$, we have 

\begin{enumerate}[(a)]
\item If $p \in ||G^{\ast}_0||$, $\sigma \in {\mathcal{S}}(G)$ if and only if either $h \geq 2$ or $m_1 \geq 2$. 

\item If $||G^{\ast}_0|| = \{ p^2 \}$, $\sigma \in {\mathcal{S}}(G)$ if and only if $h \geq 2$.
\end{enumerate}

\item For $N \geq 2$, we have

\begin{enumerate}[(a)$^{\prime}$]

\item while $||G^{\ast}_0|| = \{ p \}$, $\sigma \in {\mathcal S}(G)$ if and only if either $h \geq 2$ or $m_1 \geq 2$. 

\item while $||G^{\ast}_0|| = \{ p^2 \}$, we have $\sigma \in {\mathcal S}(G)$ if and only if either $h \geq 2$ or $m_2 \geq 2$.

\item While $||G^{\ast}_0|| = \{ p, p^2 \}$, we have $\sigma \in {\mathcal S}(G)$ if and only if either of $h \geq 2$, $m_1 \geq 2$ or, $m_2 \geq 2$ holds.
\end{enumerate}
\end{enumerate}
\end{theorem}
 
\noindent {\bf Proof} (i) is immediate as $G$ is non-abelian and minimally $2$-generated. 

\noindent (ii) The necessity of (a) follows since $G$ is minimally $2$-generated, which implies $\sigma \not\in {\mathcal{S}}(G)$ while $h = 1, m_1 = 1$. We need to prove while $||G^{\ast}_0|| = \{ p^2 \}$ and $h=1$, $\sigma \not\in {\mathcal{S}}(G)$. Suppose $\sigma = (1; m_1, 0, \dotsc, 0) \in {\mathcal{S}}(G)$ realized by $\Lambda$ as in (A2.1). Since $||G^{\ast}_1|| = \{ p^e \}$, we have $E(\Lambda) \subseteq G_2$. As $G_2 = \Phi(G)$, we have $G = \langle a_1, b_1 \rangle$. This mean $[a_1, b_1] \in G^{\ast}_1$, a contradiction.

\noindent For sufficiency of (ii)(a) while $h \geq 2$, let $s \in G^{\ast}_0, s_1 \in G^{\ast}_1$ and using {\ref{exp-descend}} choose $x \in G_1$ so that $|[s, x]| = p$. Then construct the generating set $\Lambda = \{ a_1 = s_1 = b_1, a_2 = s, b_2 = x, X_{11} = [s,x]^{-1} \}$. While $p \in ||G^{\ast}_0||$ and $m_1 \geq 2$, choose any $s \in G^{\ast}_0, s_1 \in G^{\ast}_1$ with $|s| = p$. Then construct the generating set $\Lambda_1 = \{ a_1 = s_1 = b_1, X_{11} = s, X_{12} = s^{-1} \}$ and extend by (A2.2).

\smallskip

\noindent (iii) The necessity follows from {\ref{gen1}}. We need to check the sufficiency. While $h \geq 2$, choose $s \in G^{\ast}_0, s_1 \in G^{\ast}_1$ and $x \in G_1$ so that $|[s, x]| = p^N$. Then construct the generating set $\Lambda = \{ a_1 = s_1 = b_1, a_2 = s, b_2 = x, X_{N1} = [s,x]^{-1} \}$ and extend if necessary. While $h=1, N \geq 3$, with choose $s \in G^{\ast}_0, s_1 \in G^{\ast}_1$ and $x \in G_1$ with $|x| = p^N$ and construct $\Lambda_1 = \{ a_1 = s_1 = b_1, X_{i1} = s, X_{i2} = s^{-1} x, X_{N1} = x^{-1} \}$ where $i=1$ or $2$ depending on $p$ or $p^2 \in ||G^{\ast}_0||$. While $h=1, N = 2$, corresponding to $m_1 \geq 2$ we may use $\Lambda_1$ again. While $h=1, N = 2$, corresponding to $m_2 \geq 2$ we choose $s \in G^{\ast}_0, s_1 \in G^{\ast}_1$ and construct $\Lambda_3 = \{ a_1 = s_1 = b_1, X_{21} = s, X_{22} = s^{-1} \}$. Finally, extend these by (A2.2).
\QED

\begin{table}[h]
  \begin{center}
  \caption {} \label{ed1table}
\begin{tabular}{c   c   c   c}
\hline
$N$ & $h$ & $m_N$ & property of $G$ \\
\hline
$0$      & $\geq 2$     & $0$ &  -  \\
$e$      & $\geq 1$     & $\geq 2$  &  -  \\
		 & 			   & $1, m_1 \geq 2$	 &  $||G^{\ast}_0|| = \{ p \}$  \\
		 & 		       & $1, m_2 \geq 2$ &  $||G^{\ast}_0|| = \{ p^2 \}$  \\
		 & 		       & $1, m_1 + m_2 \geq 2$ &  $||G^{\ast}_0|| = \{ p, p^2 \}$, not exceptional  \\
		 & 		       & $1, m_1 + m_2 \geq 2, m_2 \geq 1$ &  $||G^{\ast}_0|| = \{ p, p^2 \}$, $p$-exceptional  \\
		 & 		       & $1, m_1 + m_2 \geq 2, m_1 \geq 1$ &  $||G^{\ast}_0|| = \{ p, p^2 \}$, $p^2$-exceptional  \\		 
$2 \leq N \leq e-1$ & $h \geq 2$	  & - &  -  \\
		 & $h=1$	  & $m_1 \geq 2$ &  $||G^{\ast}_0|| = \{ p \}$  \\	
		 & $h=1$  & $m_2 \geq 2$ &  $||G^{\ast}_0|| = \{ p^2 \}$  \\ 
		 & $h=1$  & $m_1 \geq 2$ or $m_2 \geq 2$ &  $||G^{\ast}_0|| = \{ p, p^2 \}$  \\	
$1$ 		& $h \geq 2$	  & - &  -  \\	
		& $h=1$	  & $m_1 \geq 2$ &  $p \in ||G^{\ast}_0||$  \\	       
\hline 
\end{tabular}
\end{center}
\end{table}

\begin{remark}\label{h0note} Finally we come to the part when $\sigma \in {\mathcal {S}}(G)$ has $h=0$. While $ed(G)=1$, and $|G| \geq p^{p+2}$, we have exp($G$) $\geq p^3$. Thus all results (Thm. {\ref{0-genus-m_e-one}} and {\ref{0-genus-m_e-0}}) corresponding to $h=0$ as discussed in the previous section is identical. Table (\ref{ed1table}) summarize the results above, and for $h = 0$ see table (\ref{ed2table-expp2}). 
\end{remark}

\vspace*{.1in}

\section{\bf Results on groups of order $p^{p+1}$}\label{orderp^p+1}

The groups of order $p^{p+1}$ necessarily has ${\mathrm{exp}}(G) = p^2$. This case is somewhat special as the usual periodicity of exponent doesn't hold here. In this case we have ${\mathrm{exp}}(G_2) = p$ by {\ref{p+1-lemma}}. However, both possibilities ${\mathrm{exp}}(G_1) = p$ and ${\mathrm{exp}}(G_1) = p^2$ are possible. While ${\mathrm{exp}}(G_1) = p$ we have ${\mathrm{ed}}(G) = p \geq 2$ and while ${\mathrm{exp}}(G_1) = p^2$, we have ${\mathrm{ed}}(G) = 1$.  

While $|G| = p^{p+1}$, we have to incorporate both scenarios $c(G) = 0$ and $c(G) \geq 1$. Here we denote $G^{\ast\ast}_0 := G \setminus (G_1 \cup C_G(G_{p-1}))$. The rest of the notations $G^{\ast}_i$ ($i \geq 1$) is same as before. Note that as before we have $||G^{\ast}_1|| = \{ p^2 \}$ if and only if ${\mathrm{exp}}(G_1) = p^2$. While $c(G) = 0$, we need the following observation.

\begin{lemma}\label{cG=0-conjugacy} Let $G$ be $p$-group of maximal class of order $p^{p+1}$ and $c(G)=0$ and $\epsilon \in \{ 1, 2 \}$. If $p^{\epsilon} \in ||G^{\ast}_0 \setminus G^{\ast\ast}_0||$, then $||G^{\ast}_0 \setminus G^{\ast\ast}_0|| = \{ p^{\epsilon} \}$.
\end{lemma}

\noindent {\bf Proof.} First note that the maximal subgroup $C_G(G_{p-1})$ has order $p^p$ and hence is regular. Now the normal set $G^{\ast}_0 \setminus G^{\ast\ast}_0 = C_G(G_{p-1}) \setminus G_2$. If $s^{\prime} \in C_G(G_{p-1}) \setminus G_2$ we have
\[
C_G(G_{p-1}) \setminus G_2 = \bigcup_{j=1}^{p-1} {s^{\prime}}^j G_2
\]
Now suppose $|s^{\prime}| = p^2$ for some $s^{\prime} \in C_G(G_{p-1}) \setminus G_2$. Then $|{s^{\prime}}^j| = p^2$ for every $1 \leq j \leq p-1$. Since $C_G(G_{p-1})$ is regular and ${\mathrm{exp}}(G_2) = p$, all the elements of the normal set ${s^{\prime}}^jG_2$ has order $p^2$. The proof is same while $|s^{\prime}| = p$ as all the elements of ${s^{\prime}}^jG_2$ are non-trivial. \QED

\begin{theorem}\label{ord-p+1-h1-rel} Let Let $G$ be $p$-group of maximal class of order $p^{p+1}$ and $\sigma = (h;m_1,m_2) \in {\mathbb N}^3_{\geq 0}$. 
\begin{enumerate}[(i)]

\item If $h \geq 1$ and $m_2 \geq 2$, then $\sigma \in {\mathcal{S}}(G)$.

\item If $h \geq 1$ and $m_2 = 1$, then 

\begin{enumerate}[(a)]

\item while $||G^{\ast}_0|| = \{ p^2 \}$, then $\sigma \not\in {\mathcal S}(G)$.

\item while $||G^{\ast}_0|| = \{ p, p^2 \}$, ${\mathrm{exp}}(G_1) = p^2$ and $G$ is $p$-exceptional, then $\sigma \not\in {\mathcal S}(G)$.

\item for the remaining cases, $\sigma \in {\mathcal S}(G)$ if and only if $m_1 \geq 2$.

\end{enumerate}

\item If $h \geq 1, m_2 = 0, m_1 \geq 1$, then $\sigma \in {\mathcal{S}}(G)$.

\end{enumerate}
\end{theorem}

\noindent {\bf Proof :} (i) If ${\mathrm{exp}}(G_1) = p^2$, then choose $x \in G^{\ast}_0, y \in G^{\ast}_1$ with $|y| = p^2$. Now $\Lambda = \{ a_1 = x = b_1, X_{21} = y, X_{21} = y^{-1} \}$ realize the minimal signature $(1;0,2) \in {\mathcal{S}}(G)$. If ${\mathrm{exp}}(G_1) = p$, then we must have $p^2 \in G^{\ast}_0$. Now choose $y \in G^{\ast}_0, x \in G^{\ast}_1$ with $|y| = p^2$ and again $\Lambda$ realize $(1;0,2) \in {\mathcal{S}}(G)$.

\smallskip

\noindent (ii)(a) Let $||G^{\ast}_0|| = \{ p^2 \}$, and suppose $(h;m_1,1) \in {\mathcal{S}}(G)$ with $h \geq 1$ realized by $\Lambda$. As $E_2(\Lambda) = \{ X_{21} \}$, from long relation we have $X_{21} \not\in G^{\ast}_0$. Since ${\mathrm{exp}}(G_2) = p$, this forces $X_{21} \in G^{\ast}_1$ and consequently $||G^{\ast}_1|| = \{ p^2 \}$. But then $X_{1j} \in G_2$ for every $1 \leq j \leq m_1$ and we need at least two elliptic generators from $G^{\ast}_1$, using long relation. This is a contradiction. 

\smallskip

\noindent (ii)(b) Assume the hypothesis and suppose $(h;m_1,1) \in {\mathcal{S}}(G)$ realized by the generating set $\Lambda$ (the condition $h \geq 1$ does not play much role here). We have now $||G^{\ast}_1|| = \{ p^2 \}$. 

\smallskip

\noindent First assume $X_{21} \in G^{\ast}_1$. Then some of the $X_{1j}$ must be from $G^{\ast}_0$ (since $G_1$ is regular) and the rest must belong to $G_2$. Without loss of generality assume $X_{11}, \dotsc, X_{1k}$ are the only elements from $G^{\ast}_0$ among the elliptic elements. The long relation takes the form $X_{11} \dotsc X_{1k} X_{21} \equiv 1$ mod $G_2$. Since $G$ is $p$-exceptional, there exists integers $1 \leq \lambda_j \leq p-1$ for each $j = 2, \dotsc, k$ so that $X_{1j} = X^{\lambda_j}_{11} z_j$ for some $z_j \in G_2$. Then the long relation takes the form $X^{1 + \lambda_2 + \dotsc + \lambda_k}_{11} X_{21} \in G_2$. Then we have $1 + \lambda_2 + \dotsc + \lambda_k \equiv 0$ mod $p$. But this also means $X^{1 + \lambda_2 + \dotsc + \lambda_k}_{11} \in \Phi(G) = G_2$ forcing $X_{21} \in G_2$, a contradiction.

\smallskip

\noindent Thus we have $X_{21} \not\in G^{\ast}_1$. This mean $X_{21} \in G^{\ast}_0$. Labelling the elliptic elements of order $p$ from $G^{\ast}_0$ as above we arrive at $X^{1 + \lambda_2 + \dotsc + \lambda_k}_{11} X_{21} \in G_2$. In this case if $1 + \lambda_2 + \dotsc + \lambda_k \equiv 0$ mod $p$, then as seen above $X_{21} \in G_2$, a contradiction. Thus $p \nmid 1 + \lambda_2 + \dotsc + \lambda_k$. Now if $X_{11} \in G^{\ast\ast}_0$, it would make $X^{1 + \lambda_2 + \dotsc + \lambda_k}_{11}$ and $X^{-1}_{21}$ conjugate to each other (as these would become uniform elements), which is not possible. So $X_{11} \in C_{G}(G_{p-1}) \setminus G_2$ which imply $||G^{\ast}_0 \setminus G^{\ast\ast}_0|| = \{ p \}$ and $X_{21} \in G^{\ast}_0 \setminus G^{\ast\ast}_0$, a contradiction. 

\smallskip

\noindent (ii)(c) Now assume $p \in ||G^{\ast}_0||$ and that $(h;m_1,1) \in {\mathcal{S}}(G)$ with $h \geq 1$. We need to show $m_1 \neq 0, 1$. If $m_1 = 0$, since ${\mathrm{exp}}(G_2) = p$ we have $X_{21} \in G^{\ast}_0 \cup G^{\ast}_1$, contradicting the long relation. Now suppose $m_1 = 1$. In this case using the long relation we have $X_{21} = X^{-1}_{11}z$ for some $z \in G_2$. By \ref{unif-prop}, neither of $X_{21}, X^{-1}_{11}$ are uniform elements. Since ${\mathrm{exp}}(G_2) = p$, this means $X_{21}, X_{11} \in G^{\ast}_0 \setminus G^{\ast\ast}_0$. Now if $X_{21} \in G^{\ast}_0 \setminus G^{\ast\ast}_0$, by \ref{cG=0-conjugacy}, we have $||G^{\ast}_0 \setminus G^{\ast\ast}_0|| = \{ p^2 \}$, contradicting $X_{11} \in G^{\ast}_0 \setminus G^{\ast\ast}_0$. Similarly $X_{11} \not\in G^{\ast}_0 \setminus G^{\ast\ast}_0$. This shows $m_1 = 1$ is not possible. 

\smallskip

\noindent Now we prove sufficiency in the cases $||G^{\ast}_0|| \neq \{ p^2 \}$ and if $||G^{\ast}_0|| = \{ p, p^2 \}, {\mathrm{exp}}(G_1) = p^2$, then $G$ is not $p$-exceptional. In each of these situations we choose $x \in G^{\ast}_0, y \in G^{\ast}_1$ and set $\Lambda_1 = \{ x^{-1}, (xy), y^{-1} \}$ with various choices to realize $(0;2,1) \in {\mathcal{S}}(G)$.

\smallskip

\noindent {\it Case I :} $||G^{\ast}_0|| = \{ p \}$

\smallskip

\noindent Here we must have $||G^{\ast}_1|| = \{ p^2 \}$. Then $|x| = p = |xy|, |y| = p^2$. 

\smallskip

\noindent {\it Case II :} $||G^{\ast}_0|| = \{ p, p^2 \}, {\mathrm{exp}}(G_1) = p^2$, $G$ not $p$-exceptional

\smallskip

\noindent Choices on $x, y$ can be made so that $|x| = p = |xy|$ as $G$ is not $p$-exceptional. Now $|y| = p^2$. 

\smallskip

\noindent {\it Case III :} $||G^{\ast}_0|| = \{ p, p^2 \}, {\mathrm{exp}}(G_1) = p$, $G$ is $p$-exceptional

\smallskip

\noindent Choices can be made on $x, y$ with $|x| = p$ and $|xy| = p^2$. Now $|y| = p$.

\smallskip

\noindent {\it Case IV :} $||G^{\ast}_0|| = \{ p, p^2 \}, {\mathrm{exp}}(G_1) = p$, $G$ is $p^2$-exceptional or non exceptional

\smallskip

\noindent Choices can be made with $|x| = p^2$ and $|xy| = p$. Now $|y| = p$. 

\smallskip

\noindent (iii) Choose $a_1 \in G^{\ast}_0, b_1 \in G^{\ast}_1$ and let $X_{11} := [b_1, a_1]$. Since ${\mathrm{exp}}(G_2) = p$, the generating set $\Lambda_2 = \{ a_1, b_1, X_{11} \}$ realize the minimal signature $(1;1,0)$. \QED

\begin{theorem} Let Let $G$ be $p$-group of maximal class of order $p^{p+1}$ and $\sigma = (0;m_1,m_2) \in {\mathbb N}^3_{\geq 0}$ with $m_1 + m_2 \geq 3$ and $m_2 \geq 1$. 
\begin{enumerate}[(i)]

\item If $||G^{\ast}_0|| = \{ p \}$ then $\sigma \in {\mathcal{S}}(G)$ if and only if $m_1 \geq 2$.

\item If $||G^{\ast}_0|| = \{ p^2 \}$ and  

\begin{enumerate}[(a)]

\item if $||G^{\ast}_1|| = \{ p \}$, then $\sigma \in {\mathcal{S}}(G)$ if and only if $m_2 \geq 2$.

\item if $||G^{\ast}_1|| = \{ p^2 \}$, then $\sigma \in {\mathcal{S}}(G)$ if and only if $m_2 \geq 3$.

\end{enumerate}

\item If $||G^{\ast}_0|| = \{ p, p^2 \}$ and $G$ is $p$-exceptional and

\begin{enumerate}[(a)]

\item if $||G^{\ast}_1|| = \{ p \}$ and $p \geq 5$, then $\sigma \in {\mathcal{S}}(G)$. 

\item if $||G^{\ast}_1|| = \{ 3 \}$, then $\sigma \in {\mathcal{S}}(G)$ if and only if either $m_1 \geq 1$ or $m_1 = 0$ and $m_2 \geq 4$.

\item if $||G^{\ast}_1|| = \{ p^2 \}$, then $\sigma \in {\mathcal{S}}(G)$ if and only if $m_2 \geq 2$.

\end{enumerate}

\item If $||G^{\ast}_0|| = \{ p, p^2 \}$ and $G$ is $p^2$-exceptional and 

\begin{enumerate}[(a)]

\item ${\mathrm{exp}}(G_1) = p$, then $\sigma \in {\mathcal{S}}(G)$ if and only if $m_1 \geq 2$. 

\item ${\mathrm{exp}}(G_1) = p^2$, then $\sigma \in {\mathcal{S}}(G)$ if and only if either $m_1 \geq 1$ or $m_1 = 0$ and $m_2 \geq 4$.

\end{enumerate}

\item If $||G^{\ast}_0|| = \{ p, p^2 \}$ and $G$ is not exceptional and

\begin{enumerate}[(a)]

\item if ${\mathrm{exp}}(G_1) = p$, then $\sigma \in {\mathcal{S}}(G)$ if $G$ is of $(p^2,p^2,p^2)$-type; otherwise $\sigma \in {\mathcal{S}}(G)$ if $m_1 \geq 1$ or while $m_1 = 0$ and $m_2 \geq 4$.

\item if ${\mathrm{exp}}(G_1) = p^2$, then $\sigma \in {\mathcal{S}}(G)$

\end{enumerate}

\end{enumerate}

\end{theorem}

\noindent {\bf Proof :} (i) We have $||G^{\ast}_1|| = \{ p^2 \}$. As $\Lambda \not\subseteq G^{\ast}_1$ we have $(0;0,m_2) \not\in {\mathcal{S}}(G)$. Next if $(0;1,m_2) \in {\mathcal{S}}(G)$, then $|E(\Lambda) \cap G^{\ast}_0| = 1$ and the rest would belong to $G^{\ast}_1$ which contradict the long relation. Hence $m_1 \geq 2$ is necessary. Now the minimal signature $(0;2,1) \in {\mathcal{S}}(G)$ as seen in \ref{ord-p+1-h1-rel}.

\smallskip

\noindent (ii) If $(0;m_1,1) \in {\mathcal{S}}(G)$, then using long relation $X_{21} \not\in G^{\ast}_0$. But this means $\Lambda \subseteq G_1$, a contradiction. Now if in addition, $||G^{\ast}_1|| = \{ p^2 \}$ and $(0;m_1,2) \in {\mathcal{S}}(G)$, then $X_{11}, \dotsc, X_{1,m_1} \in G_2$ and hence $X_{22} = X^{-1}_{21} z$ for some $z \in G_2$. This contradicts $G$ is minimally $2$-generated. Now choose $x \in G^{\ast}_0, y \in G^{\ast}_1$ and consider $\Lambda_1 = \{ x^{-1}, (xy), y^{-1} \}$. This realize $(0;1,2)$ in (a) and $(0;0,3)$ in (b).

\smallskip

\noindent (iii)(a) Here $(0;2,1)$ is realized by choosing $x \in G^{\ast}_0, y \in G^{\ast}_1$ with $|x| = p = |y|$ and $|xy| = p^2$ and setting $\Lambda_1 = \{ x^{-1}, (xy), y^{-1} \}$. For $(0;1,2)$ the same generating set $\Lambda_1$ with conditions $|x| = p^2, |y| = p, |xy| = p^2$ would work. For $p \geq 5$ and $(0;0,3)$ choose $x \in G^{\ast}_0, y \in G^{\ast}_1$ with $|x| = p$. Then $|xy| = p^2 = |xy^2|$ and $(xy)(xy^2) \equiv x^2 y^3$ mod $G_2$, which implies the product represent a $z$-class in $G^{\ast}_0$ not containing $x$. Then $\Lambda_2 = \{ xy, xy^2, (xyxy^2)^{-1} \}$ realize $(0;0,3)$.   

\smallskip

\noindent (iii)(b) Here $||G^{\ast}_0|| = \{ 3, 3^2 \}$, $G$ is $3$-exceptional and $G^{\ast}_1 = \{ 3 \}$. Then $(0;2,1), (0;1,2) \in {\mathcal{S}}(G)$ as above. While $m_1 = 0$, $\Lambda = E_2(\Lambda) \subseteq G^{\ast}_0$. Now if $(0;0,3) \in {\mathcal{S}}(G)$, using the long relation if two of the generators belong to the same $z$-class in $G^{\ast}_0$, then so would be the third. Then $\Lambda$ is contained in a maximal subgroup, a contradiction. So these generators belong to distinct $z$-classes in $G^{\ast}_0$. But $G$ has only three $z$-classes and exactly one of them contain elements of order $3$, a contradiction. Now the minimal signature $(0;0,4)$ is achieved by choosing $x \in G^{\ast}_0, y \in G^{\ast}_1$ with $|x| = 3$ and hence $|xy| = 9 = (xy^2)$ and considering the generating set $\Lambda_3 = \{ xy, (xy^2), (xy^2)^{-1}, (xy)^{-1} \}$.

\smallskip

\noindent (iii)(c) As seen in \ref{ord-p+1-h1-rel}, $m_2 \geq 2$ is necessary. Now the minimal signatures $(0;1,2)$ and $(0;0,3)$ are realized by choosing $x \in G^{\ast}_0, y \in G^{\ast}_1$ with conditions $(|x|, |y|, |xy|) = (p,p^2,p^2)$ and $(p^2,p^2,p^2)$ respectively and setting the generating set $\Lambda_4 = \{ x^{-1}, (xy), y^{-1} \}$.

\smallskip

\noindent (iv)(a) As $G$ is $p^2$-exceptional, $G^{\ast}_0$ has only one $z$-class $C$ with order $p^2$ elements. If $(0;0,m_2) \in {\mathcal{S}}(G)$, then $E_2(\Lambda) \subseteq C$ and hence contained in a maximal subgroup, a contradiction. We need to check $(0;1,m_2) \not\in {\mathcal{S}}(G)$: if this is not true, then by hypothesis $X_{2j} = X_{21}^{\lambda_j} z_j$ for some $1 \leq \lambda_j \leq p-1$ and $z_j \in G_2$ for every $2 \leq j \leq m_2$. Using long relation we have $X_{11} X_{21}^{1 + \lambda_2 + \dotsc + \lambda_{m_2}} \in G_2$. Now if $1 + \lambda_2 + \dotsc + \lambda_{m_2} \equiv 0$ mod $p$, then $X_{21}^{1 + \lambda_2 + \dotsc + \lambda_{m_2}} \in \Phi(G) = G_2$ and hence $X_{11} \in G_2$. But as argued above elements of order $p^2$ cannot generate $G$. This implies $p \nmid 1 + \lambda_2 + \dotsc + \lambda_{m_2}$. Next, if $X_{21} \in G^{\ast\ast}_0$, then $X_{21}^{1 + \lambda_2 + \dotsc + \lambda_{m_2}}$ is an uniform element which is conjugate to $X^{-1}_{11}$ by \ref{unif-prop}, a contradiction. Thus $X_{21}^{1 + \lambda_2 + \dotsc + \lambda_{m_2}} \in G^{\ast}_0 \setminus G^{\ast\ast}_0$ and hence $||G^{\ast}_0 \setminus G^{\ast\ast}_0|| = \{ p^2 \}$ by \ref{cG=0-conjugacy}. But this also means $X_{11} \in G^{\ast}_0 \setminus G^{\ast\ast}_0$, a contradiction. Hence we must have $m_1 \geq 2$. Finally the minimal signature $(0;2,1) \in {\mathcal{S}}(G)$ is achieved by choosing $x \in G^{\ast}_0, y \in G^{\ast}_1$ with $|x| = p^2, |y| = p, |xy| = p$ and setting $\Lambda_4$ as above. 

\smallskip

\noindent (iv)(b) We first need to show that $(0;0,3) \not\in {\mathcal{S}}(G)$: if this is not true, then $X_{21} X_{22} X_{23} = 1$. If $E_2(\Lambda) \subseteq G^{\ast}_0$, then using $p^2$-exceptional property, they would belong to a maximal subgroup. Thus $|E_2(\Lambda) \cap G^{\ast}_0| = 2, |E_2(\Lambda) \cap G^{\ast}_1| = 1$, and this contradict $G$ is $p^2$-exceptional. Now we need to realize the signatures $(0;1,2), (0;2,1), (0;0,4)$: choose $x \in G^{\ast}_0, y \in G^{\ast}_1$ with $|x| = p^2 = |y|$. Then $|xy| = p$ and setting $\Lambda_4$ as above achieve $(0;1,2) \in {\mathcal{S}}(G)$. Next, with the same choices of $x, y$ we have $|xy| = p = |xy^{p-1}|$ using $p^2$-exceptional property. But $(xy)(xy^{p-1}) \equiv x^2y^p \equiv x^2$ mod $G_2$ and $x^2$ belong to the $z$-class that contain $x$. Hence $\Lambda_5 = \{ (xy), (xy^{p-1}), \left( (xy)(xy^{p-1}) \right)^{-1} \}$ realize $(0;2,1) \in {\mathcal{S}}(G)$. Finally using the same choices of $x, y$, $\Lambda_6 = \{ x, x^{-1}, y, y^{-1} \}$ realize $(0;0,4) \in {\mathcal{S}}(G)$. 

\smallskip

\noindent (v)(a) First choose $x \in G^{\ast}_0, y \in G^{\ast}_1$ with conditions $(|x|, |y|, |xy|) = (p^2,p,p^2), (p,p,p^2)$ respectively. Then $\Lambda_4$ as above realize $(0;1,2)$ and $(0;2,1)$ respectively. Now if $G$ is not exceptional and of $(p^2,p^2,p^2)$-type, then the signature $(0;0,3)$ is realized with all three generators belong to distinct $z$-classes in $G^{\ast}_0$ containing elements of order $p^2$. So, if $G$ is not of $(p^2,p^2,p^2)$-type, then choose $x \in G^{\ast}_0, y \in G^{\ast}_1$ with $|x| = p^2 = |xy|$ and $\Lambda_7 = \{ x, xy, (xy)^{-1}, x^{-1} \}$ realize $(0;0,4)$.    

\smallskip

\noindent (v)(b) All three minimal signatures $(0;2,1), (0;1,2), (0;0,3)$ are realized by choosing $x \in G^{\ast}_0, y \in G^{\ast}_1$ with choices $(|x|, |y|, |xy|) = (p,p^2,p), (p,p^2,p^2), (p^2,p^2,p^2)$ respectively and setting $\Lambda_4$ as above. \QED

\begin{theorem} Let Let $G$ be $p$-group of maximal class of order $p^{p+1}$ and $\sigma = (0;m_1,0) \in {\mathbb N}^3_{\geq 0}$ with $m_1 \geq 3$. 
\begin{enumerate}[(i)]

\item If $||G^{\ast}_0|| = \{ p \}$, then $\sigma \in {\mathcal{S}}(G)$.

\item If $||G^{\ast}_0|| = \{ p^2 \}$, then $\sigma \not\in {\mathcal{S}}(G)$.

\item If $||G^{\ast}_0|| = \{ p, p^2 \}$, $G$ is $p$-exceptional, and 

\begin{enumerate}[(a)]

\item if ${\mathrm{exp}}(G_1) = p$, then $\sigma \in {\mathcal{S}}(G)$ if and only if $m_1 \geq 4$.

\item if ${\mathrm{exp}}(G_1) = p^2$, then $\sigma \not\in {\mathcal{S}}(G)$.

\end{enumerate}

\item If $||G^{\ast}_0|| = \{ p, p^2 \}$, $G$ is $p^2$-exceptional, and 

\begin{enumerate}[(a)]

\item if ${\mathrm{exp}}(G_1) = p$, then $\sigma \in {\mathcal{S}}(G)$.

\item if $p \geq 5$, ${\mathrm{exp}}(G_1) = p^2$, then $\sigma \in {\mathcal{S}}(G)$. 

\item if $p = 3$, ${\mathrm{exp}}(G_1) = 3^2$, then $\sigma \in {\mathcal{S}}(G)$ if and only if $m_1 \geq 4$.

\end{enumerate}

\item If $||G^{\ast}_0|| = \{ p, p^2 \}$, $G$ is not exceptional, and 

\begin{enumerate}[(a)]

\item if ${\mathrm{exp}}(G_1) = p$, then $\sigma \in {\mathcal{S}}(G)$.

\item if ${\mathrm{exp}}(G_1) = p^2$, then $\sigma \in {\mathcal{S}}(G)$ if and only if $G$ is of $(p,p,p)$-type; otherwise $\sigma \in {\mathcal{S}}(G)$ if and only if $m_1 \geq 4$.

\end{enumerate}

\end{enumerate}
\end{theorem}

\begin{table}[h]
  \begin{center}
  \caption {} \label{ordG=p+1-table}
  \small
\begin{tabular}{c   c   c   c}
\hline
$h$ & $N$ & $m_N$ & property of $G$ \\
\hline
$\geq 2$      & $0$     & $0$ &  -  \\
$1$ &      $2$  &  $\geq 2$  & -  \\
		 &           &  $1, m_1 \geq 2$  & $||G^{\ast}_0|| \neq \{ p^2 \}$, or  \\
		 &           &       & $||G^{\ast}_0|| = \{ p, p^2 \}$, ${\mathrm{exp}}(G_1) = p^2$, $G$ not $p$-exceptional  \\
		 &      $1$  &  -  & -  \\
$0$	     &     $2$   & $m_1 \geq 2, m_2 \geq 1$ & $||G^{\ast}_0|| = \{ p \}$ or, $||G^{\ast}_0|| = \{ p, p^2 \}$, $p^2$-except., ${\mathrm{exp}}(G_1) = p$  \\
		 &           & $m_1 + m_2 \geq 3, m_2 \geq 2$  &  $||G^{\ast}_0|| = \{ p^2 \}, {\mathrm{exp}}(G_1) = p$  \\
		 &           & $m_2 \geq 3$  &  $||G^{\ast}_0|| = \{ p^2 \}, {\mathrm{exp}}(G_1) = p^2$  \\ 
		 &           & $m_1 + m_2 \geq 3$ & $||G^{\ast}_0|| = \{ p, p^2 \}$, $p$-exceptional, ${\mathrm{exp}}(G_1) = p$, $p \geq 5$ \\
		 &           & $m_1 \geq 1, m_1 + m_2 \geq 3$ or $m_1 = 0, m_2 \geq 4$ & $||G^{\ast}_0|| = \{ 3, 3^2 \}$, $3$-exceptional, ${\mathrm{exp}}(G_1) = 3$ \\ 
		 &           & $m_1 + m_2 \geq 3, m_2 \geq 2$ & $||G^{\ast}_0|| = \{ p, p^2 \}$, $p$-exceptional, ${\mathrm{exp}}(G_1) = p^2$ \\  
		 &			 & $m_1 \geq 1, m_1 + m_2 \geq 3$ or $m_1 = 0, m_2 \geq 4$  &  $||G^{\ast}_0|| = \{ p, p^2 \}$, $p^2$-exceptional, ${\mathrm{exp}}(G_1) = p^2$ \\ 
		 &			 & $m_1 + m_2 \geq 3$  &  $||G^{\ast}_0|| = \{ p, p^2 \}$, not exceptional, \\
		 &			 &					   & ${\mathrm{exp}}(G_1) = p$, $(p^2,p^2,p^2)$-type  \\ 
		 &			 & $m_1 \geq 1$ or $m_1 = 0, m_2 \geq 4$  &  $||G^{\ast}_0|| = \{ p, p^2 \}$, not exceptional, \\
		 &			 &					   & ${\mathrm{exp}}(G_1) = p$, not $(p^2,p^2,p^2)$-type  \\ 
		 &           & $m_1 + m_2 \geq 3$ & $||G^{\ast}_0|| = \{ p, p^2 \}$, not exceptional, ${\mathrm{exp}}(G_1) = p^2$ \\ 
$0$		 &		$1$	&  $m_1 \geq 3$ & $||G^{\ast}_0|| = \{ p \}$  \\ 
    		 &       &    &  $||G^{\ast}_0|| = \{ p, p^2 \}$, $p^2$-exceptional, ${\mathrm{exp}}(G_1) = p$ \\ 
    		 &       &    &  $||G^{\ast}_0|| = \{ p, p^2 \}$, $p^2$-exceptional, ${\mathrm{exp}}(G_1) = p^2$, $p \geq 5$ \\ 
    		 &       &    &  $||G^{\ast}_0|| = \{ p, p^2 \}$, not exceptional, ${\mathrm{exp}}(G_1) = p$ \\ 
    		 &       &    &  $||G^{\ast}_0|| = \{ p, p^2 \}$, not except., ${\mathrm{exp}}(G_1) = p^2$, $(p,p,p)$-type \\  
    		 & 		 & $m_1 \geq 4$ & $||G^{\ast}_0|| = \{ p, p^2 \}$, $p$-exceptional, ${\mathrm{exp}}(G_1) = p$ \\
    		 &       &    &  $||G^{\ast}_0|| = \{ 3, 3^2 \}$, $3^2$-exceptional, ${\mathrm{exp}}(G_1) = 3^2$ \\
    		 &		 &	  & $||G^{\ast}_0|| = \{ p, p^2 \}$, not exceptional, ${\mathrm{exp}}(G_1) = p^2$, not $(p,p,p)$-type \\ \hline
\end{tabular}
\end{center}
\end{table} 

\noindent {\bf Proof :} (i) Choose $x \in G^{\ast}_0, y \in G^{\ast}_1$ and set $\Lambda_1 = \{ x, xy, (x^2y)^{-1} \}$. This realize the minimal signature $(0;3,0) \in {\mathcal{S}}(G)$.

\smallskip

\noindent (ii) Every order $p$ element must be from $G_1$, which is not enough to generate $G$.

\smallskip

\noindent (iii)(a) If $(0;3,0) \in {\mathcal{S}}(G)$, then $|\Lambda \cap G^{\ast}_0| \geq 2$. If $\Lambda \subseteq G^{\ast}_0$, then the generators belong to distinct $z$-classes, contradicting $G$ is $p$-exceptional. Without loss of generality assume $X_{11}, X_{12} \in G^{\ast}_0$ and $X_{13} \in G^{\ast}_1$. Since $X_{11} X_{12} X_{13} = 1$, we have $|X_{11}| = p = |X_{12} X_{13}|$, contradicting the $p$-exceptional property again. Finally the minimal signature $(0;4,0)$ is realized by choosing $x \in G^{\ast}_0, y \in G^{\ast}_1$ with $|x| = p$ and setting $\Lambda_2 = \{ x, x^{-1}, y, y^{-1} \}$.  

\smallskip

\noindent (iii)(b) If $(0;m_1,0) \in {\mathcal{S}}(G)$, then all generators must be either from the single $z$-class containing order $p$ elements, or from $G_2$, which cannot generate $G$.

\smallskip

\noindent (iv)(a) If ${\mathrm{exp}}(G_1) = p$, then choose $x \in G^{\ast}_0, y \in G^{\ast}_1$ with $|x| = p = |xy|$. The generating set $\Lambda_3 = \{ x^{-1}, xy, y^{-1} \}$ realize the minimal signature $(0;3,0)$. 

\smallskip

\noindent (iv)(b) If $p \geq 5$, choose $x \in G^{\ast}_0, y \in G^{\ast}_1$ with $|x| = p^2$. Then $|xy| = p = |xy^2|$. Now $(xy)(xy^2) \equiv x^2y^3$ mod $G_2$. Since $p \geq 5$, $x^2y^3$ belong to a $z$-class other than the one containing $x$, and hence $|(xy)(xy^2)| = p$. Now, the generating set $\Lambda_4 = \{ xy, xy^2, \left( (xy)(xy^2) \right)^{-1} \}$ realize the minimal signature $(0;3,0) \in {\mathcal{S}}(G)$. If $p = 3$, choose $x \in G^{\ast}_0, y \in G^{\ast}_1$ with $|x| = 3 = |xy|$. The minimal signature $(0;4,0)$ is realized by $\Lambda_5 = \{ x, x^{-1}, xy, (xy)^{-1} \}$. 

\smallskip

\noindent (v)(a) The minimal signature $(0;3,0)$ is achieved by choosing $x \in G^{\ast}_0, y \in G^{\ast}_1$ with $|x| = p = |xy|$ and setting $\Lambda_3$ as above. 

\smallskip

\noindent (v)(b) If $G$ is not of $(p,p,p)$-type choose $x \in G^{\ast}_0, y \in G^{\ast}_1$ with $|x| = p = |xy|$. The minimal signature $(0;4,0)$ is realized by $\Lambda_4 = \{ x, x^{-1}, xy, (xy)^{-1} \}$. \QED

\smallskip

\noindent All results for groups with $|G| = p^{p+1}$ are summarized in table (\ref{ordG=p+1-table}).  

\vspace*{.1in}

\section{\bf Computing spectrums}\label{spectrums}

\noindent To compute the genus spectrum we consider the diophantine equation given by
\begin{equation}\label{aux-eqn}
N = p^e h + \sum_{i=1}^{e} {\frac {1}{2}} (p^e - p^{e-i}) x_i 
\end{equation}

\noindent Consider the set $\Omega_e(p)$ of all solutions of $(1)$ with integer variables $h \geq 0, x_i \geq 0$ for each $i$. For such a realizable $N$ we may have the unique $p$-adic expansion  
\[
2N = a_0 + a_1 p + \dotsc + a_{e-1} p^{e-1} + a_e p^e
\]
with the coefficients satisfying $0 \leq a_i < p$ for $0 \leq i \leq e-1$ and $a_e \geq 0$. In this expression denote by $\tau(N) = j$ so that $a_j$ is the first non-zero coefficient. Now write $S_e(2N) = \sum_{i=0}^{e} a_i$ and we have the following result:

\begin{theorem} \cite[Thm.3.1]{kma} 
\[
\Omega_e(p) = \{ N \in {\mathbb N} ~:~ S_e(2N) \geq (e- \tau(N))(p-1) \}
\] 
\end{theorem}

\noindent Let $\sigma_e(p)$ denote the smallest stable solution of $(1)$, i.e., it is minimal with the condition that every $N \geq \sigma_e(p)$ satisfy $N \in \Omega_e(p)$. Then we know that:

\begin{theorem}\label{stable-soln} \cite[Cor.3.2]{kma}
\[
\sigma_e(p) = {\frac {1}{2}} \big[ e(p-1)p^e - 3(p^e-1) \big]
\]
\end{theorem}

\noindent In fact the integer $\sigma_e(p)$ has a solution in ({\ref{aux-eqn}}) given by $h = 0, x_1 = p-1, \dotsc, x_{e-1} = p-1, x_e = p-3$.  

\smallskip

\noindent {\bf Proof of Theorem A.} Let $G$ be a $p$-group of maximal class of order $p^n$ and ${\mathrm{exp}}(G) = p^e \geq p^2$, where $p$ is an odd prime. Let ${\tilde{g}} \geq 1$. Then ${\tilde{g}}$ is a reduced genus of $G$ if and only if it satisfies
\begin{equation}
{\tilde{g}} + p^e = h p^e + \sum_{i=1}^{e} {\frac {1}{2}} (p^e - p^{e-i}) m_i
\end{equation}
and the data $\sigma = (h; m_1, m_2, \dotsc, m_e)$ satisfies: 

\smallskip

\noindent (i) for ${\mathrm{ed}}(G) \geq 2$ one of the conditions of the tables (\ref{ed2table}), (\ref{ed2table-expp2}), \\
(ii) for ${\mathrm{ed}}(G)=1$ and $|G| \geq p^{p+2}$ one of the conditions of the table (\ref{ed1table}) while $h \geq 1$, or one of the conditions of the table (\ref{ed2table}) while $h=0$, \\
(iii) for $|G| = p^{p+1}$ one of the conditions of the table (\ref{ordG=p+1-table}). \\
The proof is now routine computations and follows from the above diophantine equations in the respective cases. \QED 

\smallskip

We will now compute the reduced minimum genus of the groups. We will notice in the proofs that the condition on $ed(G)$ does not play much role here.    

\begin{theorem}\label{min-gen} Let $G$ be a finite $p$-group of maximal class of order $p^n \geq p^{p+2}$ and exponent exp$(G) = p^e \geq p^3$, where $p$ is an odd prime. Then the reduced minimum genus of $G$ is given by
\[
{\tilde{\mu}}(G) = 
\begin{cases}
{\frac {1}{2}}(p-3)p^{e-1} &\mbox{if } p \in ||G^{\ast}_0||, G ~{\mathrm{not~exceptional}}, p \geq 5, \\
{\frac {1}{2}}(p^2-3)p^{e-2} &\mbox{if } ||G^{\ast}_0|| = \{ p^2 \}, \\
{\frac {1}{2}}(p^2-p-2)p^{e-2} &\mbox{if } G ~{\mathrm{is}} ~p{\textrm{-exceptional}}, \\
{\frac {1}{2}}(p^2-2p-1)p^{e-2} &\mbox{if } G ~{\mathrm{is}} ~p^2{\textrm{-exceptional}}. \\
\end{cases}
\]
\end{theorem}

\noindent {\bf Proof.} In all the cases the minimum reduced genus is represented by a signature $\sigma \in {\mathcal{S}}(G)$ with $h=0$. As noted in {\ref{h0note}}, this situation is only dependent on exp($G$) $\geq p^3$ and not on the $ed(G)$, we simply need to specify the signatures in each such scenario. We mention these signatures $\sigma_m$ here. The calculations are straightforward:
\[
\sigma_m = 
\begin{cases}
(0; 3, 0, \dotsc, 0) &\mbox{if } p \in ||G^{\ast}_0||, G ~{\mathrm{not~exceptional}}, p \geq 5, \\
(0; 0, 3, 0, \dotsc, 0) &\mbox{if } ||G^{\ast}_0|| = \{ p^2 \}, \\
(0; 1, 2, 0, \dotsc, 0) &\mbox{if } G ~{\mathrm{is}} ~p{\textrm{-exceptional}}, \\
(0; 2, 1, 0, \dotsc, 0) &\mbox{if } G ~{\mathrm{is}} ~p^2{\textrm{-exceptional}}
\end{cases}
\]
\QED

\begin{theorem}\label{min-gen} Let $G$ be a finite $3$-group of maximal class of order $3^n \geq 3^{5}$ and exponent exp$(G) = 3^e \geq 27$. Then the reduced minimum genus of $G$ is given by
\[
{\tilde{\mu}}(G) =
\begin{cases}
3^{e-1} &\mbox{if } ||G^{\ast}_0|| = \{ 3 \}, \\
3^{e-2} &\mbox{if } ||G^{\ast}_0|| = \{ 3, 9 \}, G ~{\mathrm{not~exceptional}}. 
\end{cases}
\]
\end{theorem}

\noindent {\bf Proof.} The first case is realized by $\sigma_m = (0; 4, 0, \dotsc, 0)$. The second case is realized by $\sigma_m = (0; 2, 1, 0, \dotsc, 0)$. \QED

\smallskip

Now we proceed to compute the stable upper genus. From the description of solutions of the equations as above, the tables ({\ref{ed2table}}) and ({\ref{ed1table}}) in the previous sections shows certain minimality of the signatures whose reduced genus belong to a particular congruence class mod $p^k$, which we will discuss now.

Let $G$ be a group of order $p^n$ and exponent $p^e$. Consider a signature $\sigma = (h;m_1, \dotsc, m_N, 0, \dotsc, 0) \in {\mathcal{S}}(G)$ with $m_N \geq 1$. Extending the corresponding generating set for $\sigma$ we know that any tuple $\sigma^{\prime} = (h^{\prime};m^{\prime}_1, \dotsc, m^{\prime}_N, 0, \dotsc, 0)$ with the conditions
\[
h^{\prime} \geq h, m^{\prime}_1 \geq m_1, \dotsc, m^{\prime}_N \geq m_N
\]
satisfy $\sigma^{\prime} \in {\mathcal{S}}(G)$. Now we have
\[
{\tilde{g}}(\sigma) = (h-1)p^e + {\frac {1}{2}} \sum_{i=1}^{N} (p^e - p^{e-i})m_i \in {\widetilde{\mathrm{sp}}}(G)
\]
This reduced genus is a multiple of $p^{e-N}$ and while this happen we will call it a {\bf reduced genus of class mod $p^{e-N}$}. Dividing by $p^{e-N}$ we have 
\[
p^{-(e-N)} {\tilde{g}}(\sigma) = (h-1)p^N + {\frac {1}{2}} \sum_{i=1}^{N} (p^N - p^{N-i})m_i
\]
Now any integer $L \in \Omega_N(p)$ is expressed as
\[
L = p^N t^{\prime} + {\frac {1}{2}} \sum_{i=1}^{N} (p^N - p^{N-i})x^{\prime}_i
\]
for $t^{\prime} \geq 0, x^{\prime}_i \geq 0$. Adding above equations we have
\[
p^{-(e-N)} {\tilde{g}}(\sigma) + L = (h + t^{\prime} -1)p^N + \sum_{i=1}^{N} (p^N - p^{N-i})(m_i + x_i)
\]
Multiplying by $p^{e-N}$ we obtain
\[
{\tilde{g}}(\sigma) + p^{e-N} L = (h + t^{\prime} -1)p^e + \sum_{i=1}^{N} (p^e - p^{e-i})(m_i + x_i)
\]
This shows for every $L \in \Omega_N(p)$ we have ${\tilde{g}}(\sigma) + p^{e-N} L \in {\widetilde{\mathrm{sp}}}(G)$. We will call the signatures $\sigma^{\prime}$ that arise this way belong to the {\bf cone} in ${\mathcal{S}}(G)$ generated by $\sigma$, denoted by $C(\sigma)$, i.e., 
\[
C(\sigma) = \{ (h^{\prime};m^{\prime}_1, \dotsc, m^{\prime}_N, 0, \dotsc, 0) ~:~ h^{\prime} \geq h, m^{\prime}_1 \geq m_1, \dotsc, m^{\prime}_N \geq m_N \}
\] 

From above discussion we have the following lemma. 

\begin{lemma} Let $\sigma = (h;m_1, \dotsc, m_N, 0, \dotsc, 0) \in {\mathcal{S}}(G)$ with $m_N \geq 1$. The signature $\sigma^{\prime} \in C(\sigma)$ if and only if ${\tilde{g}}(\sigma^{\prime}) - {\tilde{g}}(\sigma) \in \Omega_N(p)$. 
\end{lemma}

At this point we deduce the following corollary to the above lemma which says every integer larger than $\sigma_e(p) - p^e$ which is a multiple of $p^{e-N}$ ($1 \leq N \leq e-1$) comes from a cone consists of reduced genus of mod $p^{e-N}$ class. 

\begin{theorem}\label{modp-class} Let $G$ be a finite $p$-group of maximal class of order $p^n \geq p^{p+2}$ and exponent exp$(G) = p^e \geq p^3$. 
\begin{enumerate}[(i)]
\item If $p \geq 5$ and $ed(G) \geq 2$, then for any integer $K \geq \sigma_e(p) - p^e$ with $p \mid K$ we have $K \in {\widetilde{\mathrm{sp}}}(G)$.

\item If $p \geq 7$ and $ed(G) = 1$, then for any integer $K \geq \sigma_e(p) - p^e$ with $p \mid K$ we have $K \in {\widetilde{\mathrm{sp}}}(G)$.
\end{enumerate}
\end{theorem}

\noindent {\bf Proof.} First consider the case $ed(G) \geq 2$. The proof is divided into several cases :

{\bf Case-I :} $||G^{\ast}_0|| = \{ p \}$  

We know that $\sigma_{2,N} = (0;3,0,\dotsc,0,1(N),0,\dotsc,0) \in {\mathcal{S}}(G)$ where $1(N)$ denote $1$ at the $N$-th co-ordinate ($2 \leq N \leq e-1$). From above discussion, it is enough to check 
\[
{\tilde{g}}(\sigma_{2,N}) + p^{e-N} \sigma_N(p) \leq \sigma_e(p) - p^e = \Bigl( {\frac {e(p-1)}{2}} - {\frac {5}{2}} \Bigr) p^e + {\frac {3}{2}}
\]
To simplify the further calculations we calculate the quantity
\begin{eqnarray*}
T_1 & := & \sigma_e(p) - p^e - p^{e-N} \sigma_N(p) \\
    &  = & \Bigl( {\frac {e(p-1)-5}{2}} p^e + {\frac {3}{2}} \Bigr) - p^{e-N} \Bigl( {\frac {N(p-1)-3}{2}} p^N + {\frac {3}{2}} \Bigr) \\
    &  = & {\frac {(e-N)(p-1) - 2}{2}}p^e - {\frac {3}{2}}p^{e-N} + {\frac {3}{2}}
\end{eqnarray*}

Then we need to check, 
\begin{eqnarray*}
{\tilde{g}}(\sigma_{2,N}) & \leq & T_1 \\
{\mathrm{i.e., }} \Bigl( p^e - {\frac {3}{2}}p^{e-1} - {\frac {1}{2}}p^{e-N} \Bigr) & \leq & {\frac {(e-N)(p-1)}{2}}p^e - p^e - {\frac {3}{2}}p^{e-N} + {\frac {3}{2}}
\end{eqnarray*}

Hence the inequality holds if and only if
\[
\Bigl( {\frac {(e-N)(p-1)}{2}} - 2 \Bigr) p^e \geq -{\frac {3}{2}}p^{e-1} + p^{e-N} + {\frac {3}{2}}
\]
i.e., if and only if
\[
{\frac {(e-N)(p-1)}{2}} \geq 2 - {\frac {3}{2}} \cdot {\frac {1}{p}} + {\frac {1}{p^N}} + {\frac {3}{2}} \cdot {\frac {1}{p^e}}
\]
For $p \geq 5$, the left side of above is an integer $\geq 2$ while $N \geq 2, e \geq 3$ implies the right side is $<2$. 

Now for the case $N=1$, we have $\sigma_{2,1} = (0;3,0,\dotsc,0) \in {\mathcal{S}}(G)$. In this case it is enough to check
\[
{\tilde{g}}(\sigma_{2,1}) + p^{e-1} \sigma_1(p) \leq \sigma_e(p) - p^e = \Bigl( {\frac {e(p-1)}{2}} - {\frac {5}{2}} \Bigr) p^e + {\frac {3}{2}}
\]
The left side above is
\[
{\tilde{g}}(\sigma_{2,1}) + p^{e-1} \sigma_1(p) = {\frac {p-3}{2}} p^e
\]
The inequality holds if and only if 
\[
\Bigl( {\frac {(e-1)(p-1)}{2}} - {\frac {3}{2}} \Bigr) p^e \geq - {\frac {3}{2}}
\]
This holds while either $e \geq 3$, or $p \geq 5$.

{\bf Case-II :} $||G^{\ast}_0|| = \{ p^2 \}$  

We know that $\sigma_{2,N} = (0;0,3,0,\dotsc,0,1(N),0,\dotsc,0) \in {\mathcal{S}}(G)$ ($3 \leq N \leq e-1$). We need to check
\begin{eqnarray*}
{\tilde{g}}(\sigma_{2,N}) + p^{e-N} \sigma_N(p) & \leq & \sigma_e(p) - p^e \\ 
{\mathrm{i.e.,}}~~~ {\tilde{g}}(\sigma_{2,N}) = p^e - {\frac {3}{2}}p^{e-2} - {\frac {1}{2}}p^{e-N} & \leq & T_1 = {\frac {(e-N)(p-1)}{2}}p^e - p^e - {\frac {3}{2}}p^{e-N} + {\frac {3}{2}}
\end{eqnarray*}
This inequality holds if and only if 
\[
{\frac {(e-N)(p-1)}{2}} \geq 2 - {\frac {3}{2p^2}} + {\frac {1}{p^N}} - {\frac {3}{2p^e}}
\]
For $p \geq 5$, the left side of above is an integer $\geq 1$ while $N \geq 3, e \geq 3$ implies the right side is $<2$.

Now, for $N=2$, we consider $\sigma_{2,2} = (0;0,3,0,\dotsc,0) \in {\mathcal{S}}(G)$. We need to check
\[
{\tilde{g}}(\sigma_{2,2}) + p^{e-2} \sigma_2(p) \leq \sigma_e(p) - p^e = \Bigl( {\frac {e(p-1)}{2}} - {\frac {5}{2}} \Bigr) p^e + {\frac {3}{2}}
\]
The left side above is
\[
{\tilde{g}}(\sigma_{2,2}) + p^{e-2} \sigma_2(p) = (p-2) p^e
\]
The inequality holds if and only if 
\[
\Bigl( {\frac {(e-2)(p-1)}{2}} - {\frac {3}{2}} \Bigr) p^e \geq - {\frac {3}{2}}
\]
which is true while $e \geq 3, p \geq 5$. 

Finally for $N=1$, we consider $\sigma_{2,1} = (2;1,0,\dotsc,0) \in {\mathcal{S}}(G)$. We need to check
\[
{\tilde{g}}(\sigma_{2,1}) + p^{e-1} \sigma_1(p) \leq \sigma_e(p) - p^e = \Bigl( {\frac {e(p-1)}{2}} - {\frac {5}{2}} \Bigr) p^e + {\frac {3}{2}}
\]
The left side above is
\[
{\tilde{g}}(\sigma_{2,1}) + p^{e-2} \sigma_1(p) = {\frac {p-1}{2}} p^e + p^{e-1}
\]
The inequality holds if and only if 
\[
{\frac {(e-1)(p-1)}{2}} \geq {\frac {5}{2}} + {\frac {1}{p}} - {\frac {3}{2p^e}}
\]
The left side is an integer $\geq 4$ while $e \geq 3$ and $p \geq 5$ while the right side is $<3$. 

{\bf Case-III :} $||G^{\ast}_0|| = \{ p, p^2 \}$, $G$ is not $p^2$-exceptional.

First consider $3 \leq N \leq e-1$ and $\sigma_{2,N} = (0;1,2,0,\dotsc, 0, 1(N),0, \dotsc, 0) \in {\mathcal{S}}(G)$. We need to check
\begin{eqnarray*}
{\tilde{g}}(\sigma_{2,N}) + p^{e-N} \sigma_N(p) & \leq & \sigma_e(p) - p^e \\ 
{\mathrm{i.e.,}}~~~ {\tilde{g}}(\sigma_{2,N}) = p^e - {\frac {1}{2}} p^{e-1} - p^{e-2} - {\frac {1}{2}}p^{e-N} & \leq & T_1
\end{eqnarray*}
This inequality holds if and only if 
\[
{\frac {(e-N)(p-1)}{2}} \geq 2 - {\frac {1}{2p}} - {\frac {1}{p^2}} + {\frac {1}{p^N}} - {\frac {3}{2p^e}}
\]  
The left side is an integer $\geq 2$ for $p \geq 5$, while the right side is $< 2$.

For $N=2$, we consider $\sigma_{2,2} = (0;1,2,0, \dotsc, 0) \in {\mathcal{S}}(G)$. We need to check
\[
{\tilde{g}}(\sigma_{2,2}) + p^{e-2} \sigma_2(p) \leq \sigma_e(p) - p^e 
\]
This inequality holds if and only if
\[
{\frac {(e-2)(p-1)}{2}} \geq {\frac {3}{2}} - {\frac {1}{2p}} + {\frac {1}{2p^2}} - {\frac {3}{2p^e}}
\]
which is true for $e \geq 3, p \geq 5$. For $N=1$ we can again use $\sigma_{2,1} = (2;1,0, \dotsc,0) \in {\mathcal{S}}(G)$ as in the second case.

{\bf Case-IV :} $||G^{\ast}_0|| = \{ p, p^2 \}$, $G$ is $p^2$-exceptional.

First consider $3 \leq N \leq e-1$ and $\sigma_{2,N} = (0;2,1,0,\dotsc, 0, 1(N),0, \dotsc, 0) \in {\mathcal{S}}(G)$. We need to check
\begin{eqnarray*}
{\tilde{g}}(\sigma_{2,N}) + p^{e-N} \sigma_N(p) & \leq & \sigma_e(p) - p^e \\ 
{\mathrm{i.e.,}}~~~ {\tilde{g}}(\sigma_{2,N}) = p^e - p^{e-1} - {\frac {1}{2}} p^{e-2} - {\frac {1}{2}}p^{e-N} & \leq & T_1
\end{eqnarray*}
This inequality holds if and only if 
\[
{\frac {(e-N)(p-1)}{2}} \geq 2 - {\frac {1}{p}} - {\frac {1}{2p^2}} + {\frac {1}{p^N}} - {\frac {3}{2p^e}}
\]  
The left side is an integer $\geq 2$ for $p \geq 5$, while the right side is $< 2$. 

For $N=2$, we consider $\sigma_{2,2} = (0;2,1,0, \dotsc, 0) \in {\mathcal{S}}(G)$. We need to check
\[
{\tilde{g}}(\sigma_{2,2}) + p^{e-2} \sigma_2(p) \leq \sigma_e(p) - p^e 
\]
This inequality holds if and only if
\[
{\frac {(e-2)(p-1)}{2}} \geq {\frac {3}{2}} - {\frac {1}{p}} + {\frac {1}{p^2}} - {\frac {3}{2p^e}}
\]
which is true for $e \geq 3, p \geq 5$. For $N=1$ we can again use $\sigma_{2,1} = (2;1,0, \dotsc,0) \in {\mathcal{S}}(G)$ as in the second case. This completes the case $ed(G) \geq 2$.

Now we consider the case $ed(G) = 1$.

In this we first note that $\sigma_{2,N} = (2;0,\dotsc,0,1(N),0,\dotsc,0) \in {\mathcal{S}}(G)$ for $1 \leq N \leq e-1$ irrespective of the type of $G$. We need to check
\begin{eqnarray*}
{\tilde{g}}(\sigma_{2,N}) + p^{e-N} \sigma_N(p) & \leq & \sigma_e(p) - p^e \\ 
{\mathrm{i.e.,}}~~~ {\tilde{g}}(\sigma_{2,N}) = {\frac {3}{2}} p^e  - {\frac {1}{2}}p^{e-N} & \leq & T_1 = {\frac {(e-N)(p-1)}{2}}p^e - p^e - {\frac {3}{2}}p^{e-N} + {\frac {3}{2}}
\end{eqnarray*}
This inequality holds if and only if 
\[
{\frac {(e-N)(p-1)}{2}} \geq {\frac {5}{2}} + {\frac {1}{p^N}} - {\frac {3}{2p^e}}
\]  
While $p \geq 7$, the left side is an integer $\geq 6$ and the right side is $< {\frac {5}{2}}$. 
\QED

\begin{definition} A signature $\sigma \in {\mathcal{S}}(G)$ is called {\bf minimal} if the cone $C(\sigma) \subseteq C(\sigma^{\prime})$ implies $\sigma = \sigma^{\prime}$.
\end{definition}

\begin{theorem} Let $p \geq 5$ be an odd prime and $G$ be a finite $p$-group of maximal class of order $p^n \geq p^{p+2}$ and exponent exp$(G) = p^e \geq p^3$. Let $ed(G) \geq 2, p \in ||G^{\ast}_0||$ and $G$ is not $p$-exceptional. Then the reduced stable upper genus of $G$ is given by ${\tilde{\sigma}}(G) =  \sigma_e(p) - p^e$
\end{theorem} 

\noindent {\bf Proof.} We know from ({\ref{ed2table}}) that the (minimal) signatures $\sigma_1, \sigma_2 \in {\mathcal{S}}(G)$ given by
\[
\sigma_1 = (1;0,\dotsc,0,1), \sigma_2 = (0;2,0,\dotsc,0,1)
\] 

By {\ref{stable-soln}} every integer $N \geq \sigma_e(p)$ has a solution
\[
N = hp^e + \sum_{i=1}^{e} {\frac {1}{2}}(p^e - p^{e-i}) x_i
\]
for some $h, x_i \geq 0$. For each such  $N$ we have 
\begin{eqnarray*}
N + {\mathrm{\bf{\tilde{g}}}}(0;2,0, \dotsc,0,1) & = & {\mathrm{\bf{\tilde{g}}}}(h;2+x_1, x_2, \dotsc, x_{e-1}, x_e + 1) \\
 & \geq & \sigma_e(p) - p^e + {\frac {1}{2}} (p^e - p^{e-1}) \cdot 2 + {\frac {1}{2}}(p^e - 1)  
\end{eqnarray*}
and $(h; 2+x_1, x_2, \dotsc, x_{e-1}, x_e + 1) \in {\mathcal{S}}(G)$. This means
\begin{eqnarray*}
\sigma_e(p) - p^e \leq {\tilde{\sigma}}(G) & \leq & \sigma_e(p) - p^e + {\frac {1}{2}} (p^e - p^{e-1}) \cdot 2 + {\frac {1}{2}}(p^e - 1) \\
& = & {\frac {1}{2}} e(p-1)p^e - p^e - p^{e-1} + 1 =: M_0 
\end{eqnarray*}

We will now show that if $a$ is an integer with $a \not\equiv 0$ mod $p$ with $\sigma_e(p) - p^e + 1 \leq a \leq M_0 - 1$, then there exists a $\sigma \in C(\sigma_1) \cup C(\sigma_2)$ so that ${\tilde{g}}(\sigma) = a$. Set $a = M_0 - \alpha$ so that 
\[
1 \leq \alpha \leq M_0 - (\sigma_e(p) - p^e + 1) = {\frac {3}{2}}p^e - p^{e-1} - {\frac {3}{2}}
\]
This implies $2\alpha$ is an even integer and it can be written as a $p$-adic expansion $2\alpha = b_0 + b_1 p + \dotsc + b_{e-1} p^{e-1} + b_e p^e$ with $0 \leq b_i \leq p-1$ for $0 \leq i \leq e-1$ and $0 \leq b_e \leq 2$. Moreover, while $b_e = 2$, we have $b_{e-1} \leq p-3$.

The condition $a \not\equiv 0$ mod $p$ is equivalent to $b_0 \neq 2$. For any such $\alpha$, $M_0 - \alpha \in C(\sigma_2)$ if and only if $M_0 - \alpha - {\tilde{g}}(\sigma_2) \in \Omega_e(p)$.

We will show that if $M_0 - \alpha - {\tilde{g}}(\sigma_2) \not\in \Omega_e(p)$, then $M_0 - \alpha \in C(\sigma_1)$. Set $J := 2(M_0 - \alpha - {\tilde{g}}(\sigma_2))$. Then,
\begin{eqnarray*}
J & = & 2(M_0 - {\tilde{g}}(\sigma_2)) - 2\alpha \\
  & = & \Bigl( e(p-1) - 3 \Bigr) p^e + 3 - \Bigl( b_0 + b_1 p + \dotsc + b_{e-1} p^{e-1} + b_e p^e \Bigr) \\
  & = & \Bigl( e(p-1) - 4 - b_e \Bigr) p^e + (p-1 - b_{e-1}) p^{e-1} + \dotsc + (p-1 - b_1) p + (p+3 - b_0)
\end{eqnarray*}

First consider the case $b_e = 0,1$. Now while $4 \leq b_0 \leq p-1$ this is a $p$-adic expansion of $J$ with the coefficient of $p^e$ arbitrary. Thus, our requirement for ${\frac {1}{2}}J \in \Omega_e(p)$ is
\[
S_e(J) = 2e(p-1) - (b_0 + \dotsc + b_{e-1} + b_e) \geq e(p-1)
\] 
While $b_e = 0$ this is always true as $b_0 + \dotsc + b_{e-1} \leq e(p-1)$. While $b_e = 1$, one of $b_0, \dotsc, b_{e-1}$ must be $\leq p-2$ as $2\alpha$ is even, and hence this holds. This implies for all such $\alpha$, we have $M_0 - \alpha - {\tilde{g}}(\sigma_2) \in \Omega_e(p)$ and hence $M_0 - \alpha \in C(\sigma_2)$. 

Now consider the case $b_0 \in \{ 0, 1, 3 \}$. Let $b_k$ be the first non-zero integer among $b_1, \dotsc b_{e-1}$. Then the required $p$-adic expansion of $J$ is given by
\[
\Bigl( e(p-1) - 4 - b_e \Bigr) p^e + (p-1 - b_{e-1}) p^{e-1} + \dotsc + (p-1 - b_{k+1})p^{k+1} + (p - b_k) p^k + (3 - b_0)
\] 

Then, 
\[
S_e(J) = e(p-1) + (e-k)(p-1) - (b_0 + b_k + \dotsc + b_{e-1} + b_e)
\]
While $b_0 = 3$, we need $S_e(J) \geq (e-k)(p-1)$; equivalently, 
\[
b_k + \dotsc + b_{e-1} + b_e \leq e(p-1) - 3 = (e-1)(p-1) + (p-3) - 1
\] 
which is always true. While $b_0 = 0, 1$ we need $S_e(J) \geq e(p-1)$; equivalently,
\[
b_0 + b_k + \dotsc + b_{e-1} + b_e \leq (e-k)(p-1)
\]
In case $b_0 = 0$, either $b_e = 0$, or if $b_e = 1$, one of $b_k, \dotsc, b_{e-1}$ must be $\leq p-2$ to keep $2\alpha$ even. Hence this inequality holds. In case $b_0 = 1$, this inequality is contradicted only at the extremal cases
\[
b_0 = 1, b_k = b_{k+1} = \dotsc = b_{e-1} = p-1, b_e = 1
\]
Denote by $\alpha_k$ these values where $1 \leq k \leq e-1$. We will show that in such situation $M_0 - \alpha_k \in C(\sigma_1)$ i.e., $M_0 - \alpha_k - {\tilde{g}}(\sigma_1) \in \Omega_e(p)$. Set $J_k^{\prime} = 2(M_0 - \alpha_k - {\tilde{g}}(\sigma_1))$. While $1 \leq k \leq e-2$ we have the required $p$-adic expansion of $J_k^{\prime}$ as
\begin{eqnarray*}
J_k^{\prime} & = & 2(M_0 - {\tilde{g}}(\sigma_1)) - 2\alpha_k \\
  & = & \Bigl( e(p-1) - 3 \Bigr) p^e - 2p^{e-1} + 3 - \Bigl( 1 + (p-1) p^k + \dotsc + (p-1) p^{e-1} + p^e \Bigr) \\
  & = & \Bigl( e(p-1) - 6 \Bigr) p^e + (p-2) p^{e-1} + p^k + 2
\end{eqnarray*}

Now $S_e(J_k^{\prime}) = e(p-1) + p-5 \geq e(p-1)$, as $p \geq 5$. Similarly while $k = e-1$, we have
\[
J_{e-1}^{\prime} = \Bigl( e(p-1) - 6 \Bigr) p^e + (p-1)p^{e-1} + 2
\] 
and again $S_e(J_{e-1}^{\prime}) \geq e(p-1)$ as $p \geq 5$.

We now turn to the case $b_e = 2$. In this situation we have
\begin{eqnarray*}
J & = & 2(M_0 - {\tilde{g}}(\sigma_2)) - 2\alpha \\
  & = & \Bigl( e(p-1) - 3 \Bigr) p^e + 3 - \Bigl( b_0 + b_1 p + \dotsc + b_{e-1} p^{e-1} + b_e p^e \Bigr) \\
  & = & \Bigl( e(p-1) - 6 \Bigr) p^e + (p-1 - b_{e-1}) p^{e-1} + \dotsc + (p-1 - b_1) p + (p+3 - b_0)
\end{eqnarray*}

As before while $4 \leq b_0 \leq p-1$ this is a $p$-adic expansion of $J$ with the coefficient of $p^e$ arbitrary. Thus, our requirement for ${\frac {1}{2}}J \in \Omega_e(p)$ is
\[
S_e(J) = 2e(p-1) - 2 - (b_0 + \dotsc + b_{e-1}) \geq e(p-1)
\]
i.e., $b_0 + \dotsc + b_{e-1} \leq e(p-1) - 2$. Since $b_{e-1} \leq p-3$, this is always true. So assume $b_0 \in \{ 0, 1, 3 \}$. Let $b_k$ be the first non-zero integer among $b_1, \dotsc b_{e-1}$. Then the required $p$-adic expansion of $J$ is given by
\[
\Bigl( e(p-1) - 6 \Bigr) p^e + (p-1 - b_{e-1}) p^{e-1} + \dotsc + (p-1 - b_{k+1})p^{k+1} + (p - b_k) p^k + (3 - b_0)
\]
Then, 
\[
S_e(J) = e(p-1) - 2 + (e-k)(p-1) - (b_0 + b_k + \dotsc + b_{e-1})
\]
While $b_0 = 3$, we need $S_e(J) \geq (e-k)(p-1)$; equivalently, 
\[
b_k + \dotsc + b_{e-1} \leq e(p-1) - 5 = (e-1)(p-1) + (p-5)
\] 
which is always true. While $b_0 = 0,1$, we need $S_e(J) \geq e(p-1)$; equivalently,
\[
b_k + \dotsc + b_{e-1} \leq (e-k)(p-1) - 2 - b_0 
\]  
While $b_0 = 0$, this again follows from $b_{e-1} \leq p-3$. While $b_0 = 1$, either one of $b_1, \dotsc, b_{e-2}$ is $\leq p-2$, or $b_{e-1} \leq p-4$, since $2\alpha$ is even. Hence this follows again. 

Finally note that while $b_0 = 2$, we have $p \mid M_0 - \alpha$ and hence $M_0 - \alpha \in {\widetilde{\mathrm{sp}}}(G)$ by {\ref{modp-class}}.
\QED

\begin{theorem} Let $G$ be a finite $p$-group of maximal class of order $p^n \geq p^{p+2}$ and exponent exp$(G) = p^e \geq p^3$, where $p \geq 5$ is an odd prime. Let $ed(G) \geq 2$, $p^2 \in ||G^{\ast}_0||$ and $G$ is not $p^2$-exceptional. Then the reduced stable upper genus of $G$ is given by ${\tilde{\sigma}}(G) = \sigma_e(p) - p^e$.
\end{theorem}

\noindent {\bf Proof.} Using ({\ref{ed2table}}), as in the previous theorem we use the following two minimal signatures to analyse the reduced genus belonging to mod $1$ class :  
\[
\sigma_1 = (1;0,\dotsc,0,1), \sigma_2 = (0;0,2,0,\dotsc,0,1)
\] 

By {\ref{stable-soln}} every integer $N \geq \sigma_e(p)$ has a solution
\[
N = hp^e + \sum_{i=1}^{e} {\frac {1}{2}}(p^e - p^{e-i}) x_i
\]
for some $h, x_i \geq 0$. For each such  $N$ we have 
\begin{eqnarray*}
N + {\mathrm{\bf{\tilde{g}}}}(0;0,2,0, \dotsc,0,1) & = & {\mathrm{\bf{\tilde{g}}}}(h;x_1, 2+x_2, \dotsc, x_{e-1}, x_e + 1) \\
 & \geq & \sigma_e(p) - p^e + {\frac {1}{2}} (p^e - p^{e-2}) \cdot 2 + {\frac {1}{2}}(p^e - 1)  
\end{eqnarray*}
and $(h; x_1, 2+x_2, \dotsc, x_{e-1}, x_e + 1) \in {\mathcal{S}}(G)$. This means
\begin{eqnarray*}
\sigma_e(p) - p^e \leq {\tilde{\sigma}}(G) & \leq & \sigma_e(p) - p^e + {\frac {1}{2}} (p^e - p^{e-2}) \cdot 2 + {\frac {1}{2}}(p^e - 1) \\
& = & {\frac {1}{2}} e(p-1)p^e - p^e - p^{e-2} + 1 =: M_0 
\end{eqnarray*}

We will now show that if $a$ is an integer with $a \not\equiv 0$ mod $p$ with $\sigma_e(p) - p^e + 1 \leq a \leq M_0 - 1$, then there exists a $\sigma \in {\mathcal{S}}(G)$ so that ${\tilde{g}}(\sigma) = a$. Set $a = M_0 - \alpha$ so that 
\[
1 \leq \alpha \leq M_0 - (\sigma_e(p) - p^e + 1) = {\frac {3}{2}}p^e - p^{e-2} - {\frac {3}{2}}
\]
This implies $2\alpha$ is an even integer and it can be written as a $p$-adic expansion $2\alpha = b_0 + b_1 p + \dotsc + b_{e-1} p^{e-1} + b_e p^e$ with $0 \leq b_i \leq p-1$ for $0 \leq i \leq e-1$ and $0 \leq b_e \leq 2$. Moreover, while $b_e = 2, b_{e-1} = p-1$, we have $b_{e-2} \leq p-3$.

The condition $a \not\equiv 0$ mod $p$ is equivalent to $b_0 \neq 2$. For any such $\alpha$, $M_0 - \alpha \in C(\sigma_2)$ if and only if $M_0 - \alpha - {\tilde{g}}(\sigma_2) \in \Omega_e(p)$.

We will show that if $M_0 - \alpha - {\tilde{g}}(\sigma_2) \not\in \Omega_e(p)$, then $M_0 - \alpha \in C(\sigma_1)$. Set $J := 2(M_0 - \alpha - {\tilde{g}}(\sigma_2))$. Then,
\begin{eqnarray*}
J & = & 2(M_0 - {\tilde{g}}(\sigma_2)) - 2\alpha \\
  & = & \Bigl( e(p-1) - 3 \Bigr) p^e + 3 - \Bigl( b_0 + b_1 p + \dotsc + b_{e-1} p^{e-1} + b_e p^e \Bigr) \\
  & = & \Bigl( e(p-1) - 4 - b_e \Bigr) p^e + (p-1 - b_{e-1}) p^{e-1} + \dotsc + (p-1 - b_1) p + (p+3 - b_0)
\end{eqnarray*}

The case $b_e = 0,1$ is identical to the previous theorem. So we consider the case $b_e = 2$. In this situation we have
\begin{eqnarray*}
J & = & 2(M_0 - {\tilde{g}}(\sigma_2)) - 2\alpha \\
  & = & \Bigl( e(p-1) - 3 \Bigr) p^e + 3 - \Bigl( b_0 + b_1 p + \dotsc + b_{e-1} p^{e-1} + b_e p^e \Bigr) \\
  & = & \Bigl( e(p-1) - 6 \Bigr) p^e + (p-1 - b_{e-1}) p^{e-1} + \dotsc + (p-1 - b_1) p + (p+3 - b_0)
\end{eqnarray*}

As before while $4 \leq b_0 \leq p-1$ this is a $p$-adic expansion of $J$ with the coefficient of $p^e$ arbitrary. Thus, our requirement for ${\frac {1}{2}}J \in \Omega_e(p)$ is
\[
S_e(J) = 2e(p-1) - 2 - (b_0 + \dotsc + b_{e-1}) \geq e(p-1)
\]
i.e., $b_0 + \dotsc + b_{e-1} \leq e(p-1) - 2$. While $b_{e-1} \leq p-3$, this is always true. While $b_{e-1} = p-2$, we need another $b_i \leq p-2$ to keep $2\alpha$ even, so this holds. Finally while $b_{e-1} = p-1$ we have $b_{e-2} \leq p-3$, again the inequality holds.

So assume $b_0 \in \{ 0, 1, 3 \}$. Let $b_k$ be the first non-zero integer among $b_1, \dotsc b_{e-1}$. Then the required $p$-adic expansion of $J$ is given by
\[
\Bigl( e(p-1) - 6 \Bigr) p^e + (p-1 - b_{e-1}) p^{e-1} + \dotsc + (p-1 - b_{k+1})p^{k+1} + (p - b_k) p^k + (3 - b_0)
\]
Then, 
\[
S_e(J) = e(p-1) - 2 + (e-k)(p-1) - (b_0 + b_k + \dotsc + b_{e-1})
\]
While $b_0 = 3$, we need $S_e(J) \geq (e-k)(p-1)$; equivalently, 
\[
b_k + \dotsc + b_{e-1} \leq e(p-1) - 5 = (e-1)(p-1) + (p-5)
\] 
which is always true. While $b_0 = 0,1$, we need $S_e(J) \geq e(p-1)$; equivalently,
\[
b_k + \dotsc + b_{e-1} \leq (e-k)(p-1) - 2 - b_0 
\]  
While $b_0 = 0, b_{e-1} = p-1$, this again follows from $b_{e-2} \leq p-3$. While $b_0 = 0, b_{e-1} \leq p-2$, either $b_{e-1} \leq p-3$ or if $b_{e-1} = p-2$, one of $b_1, \dotsc, b_{e-2}$ is $\leq p-2$, to keep $2\alpha$ even.

Now consider $b_0 = 1$. While $b_0 = 1, b_{e-1} = p-1$, we have $b_{e-2} \leq p-3$. Then either one of $b_1, \dotsc, b_{e-3}$ is $\leq p-2$, or $b_{e-2} \leq p-4$ to keep $2\alpha$ even. Similarly when $b_0 = 1, b_{e-1} \leq p-3$, the inequality holds.

The only extremal case arise in this situation is
\[
b_0 = 1, b_k = \dotsc = b_{e-2} = p-1, b_{e-1} = p-2 
\]

Denote by $\alpha_k$ these values where $1 \leq k \leq e-2$. We will show that in such situation $M_0 - \alpha_k \in C(\sigma_1)$ i.e., $M_0 - \alpha_k - {\tilde{g}}(\sigma_1) \in \Omega_e(p)$. Set $J_k^{\prime} = 2(M_0 - \alpha_k - {\tilde{g}}(\sigma_1))$. Now we consider the $p$-adic expansion of $J_k^{\prime}$ as
\begin{eqnarray*}
J_k^{\prime} & = & 2(M_0 - {\tilde{g}}(\sigma_1)) - 2\alpha_k \\
  & = & \Bigl( e(p-1) - 3 \Bigr) p^e - 2p^{e-1} + 3 \\
  &   & -\Bigl( 1 + (p-1) p^k + \dotsc + (p-1)p^{e-2} + (p-2) p^{e-1} + 2p^e \Bigr) \\
  & = & \Bigl( e(p-1) - 7 \Bigr) p^e + (p-1) p^{e-1} + p^k + 2
\end{eqnarray*}
and we have
\[
S_e(J^{\prime}_k) = e(p-1) + (p-5)
\]
which is $\geq e(p-1)$. Finally, while $b_0 = 2$, we have $p \mid M_0 - \alpha$ and hence $M_0 - \alpha \in {\widetilde{\mathrm{sp}}}(G)$ by {\ref{modp-class}}.
\QED

We will now discuss the case $ed(G) = 1$. Here we require the assumption $p \geq 7$ since we have seen in {\ref{modp-class}} that some of the integers above $\sigma_e(p) - p^e$ which are $\equiv 0$ mod $p$ does not come from signatures of mod $p^k$ class. Thus either these produce gaps in the genus spectra, or are recovered by the mod $1$ class. These require some more detail combinatorics. 

\begin{theorem} Let $G$ be a finite $p$-group of maximal class of order $p^n \geq p^{p+2}$ and exponent exp$(G) = p^e \geq p^3$, where $p \geq 7$ is an odd prime. Let $ed(G)=1$. Then the reduced stable upper genus of $G$ is given by ${\tilde{\sigma}}(G) = \sigma_e(p) - p^e$
\end{theorem}

\noindent {\bf Proof.} The proof of this uses similar techniques except in all the groups we have the following two signatures which generate the cone giving almost all reduced genus of mod $1$ class 
\[
\sigma_1 = (1;0,\dotsc,0,2), \sigma_2 = (0;2,0,\dotsc,0,1)
\]
Especially, it is enough to check which integers $a$ with $\sigma_e(p)-p^e+1 \leq a \leq M_0 - 1$ are in $C(\sigma_1) \cup C(\sigma_2)$, where 
\[
M_0 = 
\begin{cases}
{\frac {1}{2}}e(p-1)p^e - p^e - p^{e-1}+1 &\mbox{if } p \in ||G^{\ast}_0||, G \neq p-{\mathrm{exceptional}} \\
{\frac {1}{2}}e(p-1)p^e - p^e - p^{e-2}+1 &\mbox{if } p^2 \in ||G^{\ast}_0||, G \neq p^2-{\mathrm{exceptional}}
\end{cases}
\]

In the first case we only need to analyse the integer $M_0 - \alpha_k$ where $2\alpha_k = b_0 + b_1 p + \dotsc + b_e p^e$ with 
\[
b_0 = 1, b_k = \dotsc = b_{e-1} = p-1, b_e = 1
\]
Setting $J^{\prime}_k = 2(M_0 - \alpha_k - {\tilde{g}}(\sigma_1))$ we have for $1 \leq k \leq e-2$ that
\begin{eqnarray*}
J^{\prime}_k & = & 2(M_0 - {\tilde{g}}(\sigma_1)) - 2\alpha_k \\
  & = & \Bigl( e(p-1) - 4 \Bigr) p^e - 2p^{e-1} + 4 - \Bigl( 1 + (p-1) p^k + \dotsc + (p-1) p^{e-1} + p^e \Bigr) \\
  & = & \Bigl( e(p-1) - 7 \Bigr) p^e + (p-2) p^{e-1} + p^k + 3
\end{eqnarray*}
Now $S_e(J_k^{\prime}) = e(p-1) + p-5 \geq e(p-1)$, as $p \geq 5$. Similarly while $k = e-1$, we have
\[
J_{e-1}^{\prime} = \Bigl( e(p-1) - 7 \Bigr) p^e + (p-1)p^{e-1} + 3
\] 
and again $S_e(J_{e-1}^{\prime}) \geq e(p-1)$ as $p \geq 7$. The second case is similar.

\begin{note} For $p=5$, we cannot indeed use the result {\ref{modp-class}} and indeed the following extremal cases arise
\[
b_0 = 2, b_k = p-1, \dotsc, b_{e-2} = p-1, b_{e-1} = p-3, b_e = 2
\]
for the first case in above theorem, which yield the stable upper genus 
\[
{\tilde{\sigma}}(G) = \Bigl( {\frac {e(p-1)-4}{2}} \Bigr) p^e - \Bigl( {\frac {p-1}{2}} \Bigr) p^{e-1} + 1
\]  
Similar calculations can be made for the other cases.
\end{note}

\section{\bf Constructing groups achieving the spectrums}\label{examples}

In this section we will compute explicit groups of given order and exponent which achieve the spectrum as seen in previous sections except the $p^2$-exceptional type. For this we go through the structure theorems of metabelian $p$-groups $G$ of maximal class (which means $G_2$ is necessarily abelian) proved by Miech. We start with the following theorem. 

\begin{theorem}\label{mie-struc} \cite[Thm.2]{mie} Let $G$ be a metabelian $p$-group of maximal class and order $p^n \geq p^{p+1}$ and $[G_1, G_2] = G_{n-k}$ for $0 \leq k \leq p-2$. Let $s \in G^{\ast}_0, s_1 \in G^{\ast}_1$ and $s_i = [s_{i-1}, s]$ for $2 \leq i \leq n-1$. Then:
\begin{enumerate}[(i)]
\item there exists integers $a_{n-k}, a_{n-k+1}, \dotsc, a_{n-1}$ with $a_{n-k} \not\equiv 0$ mod $p$ such that
\[
[s_1, s_2] = s^{a_{n-k}}_{n-k} s^{a_{n-k+1}}_{n-k+1} \dotsc s^{a_{n-1}}_{n-1}
\]  

\item there exists integers $0 \leq w, z \leq p-1$ such that
\[
s^p = s^{w}_{n-1}, s^{p \choose 1}_1 s^{p \choose 2}_2 \dotsc s^{p \choose p}_p = s^{z}_{n-1}
\] 

\item for $2 \leq i \leq n-1$ we have 
\[
s^{p \choose 1}_i s^{p \choose 2}_{i+1} \dotsc s^{p \choose p}_{i+p-1} = 1
\]
\end{enumerate}

Conversely, given a set of parameters $(k, a_{n-k}, a_{n-k+1}, \dotsc, a_{n-1}, w, z)$ there exists a metabelian $p$-group $G$ of maximal class and order $p^n$ satisfying (i)-(iii).
\end{theorem}

\begin{lemma}\label{G0lem1} \cite[Lem.8]{mie} Let $G$ be a metabelian $p$-group of maximal class of order $p^n$. Let $[G_1, G_2] = G_{n-k}$. Assuming notations of {\ref{mie-struc}} we have
\[
(s s^{\zeta}_1)^p = s^p \bigl( s^{p \choose 1}_1 s^{p \choose 2}_2 \dotsc s^{p \choose p}_p \bigr)^{\zeta} s^{\psi \zeta^2}_{n-1}
\]
where $\psi$ is given by
\[
\psi = \psi(k) = 
\begin{cases}
a_{n-k} &\mbox{if } k = p-2, \\
0 &\mbox{if } k \leq p-3
\end{cases}
\]
\end{lemma}

Given a set ${\underline {v}} = (k, a_{n-k}, a_{n-k+1}, \dotsc, a_{n-1}, w, z)$ of parameters we will consider one such group $G({\underline{v}})$ satisfying conditions of {\ref{mie-struc}} and will show that $G({\underline{v}})$ satisfy our requirements. Note that as shown in {\cite{mie}}, there could be several groups $G$ upto isomorphism which has the same parameter set ${\underline{v}}$. While $n \geq p+2$, the condition (iii) of theorem {\ref{mie-struc}} ensures that exp($G_{n-p+1}$) $=p$. We will also use the following fact which can be easily proved.

\begin{proposition}\label{abel-sub} (\cite[Exercise 4.7]{falc}) Let $G$ be a $p$-group of maximal class of order $p^n \geq p^{p+2}$. If $N \unlhd G$ is abelian of order $p^t$ with $t = k(p-1)+r, 0 \leq r \leq p-1$, then
\[
N \cong C_{p^{k+1}} \times \overset{r}{\dotsc} \times C_{p^{k+1}} \times C_{p^k} \times \overset{p-r-1}{\dotsc} \times C_{p^k} 
\] 
\end{proposition}

\begin{theorem}\label{examples-ed2} Given a prime $p \geq 3$, an integer $e \geq 3$, set $n = (e-1)(p-1)+3$. Then there exists a metabelian $p$-group $G = G({\underline{v}})$ of maximal class of order $p^n$ corresponding to the parameters ${\underline{v}}$ which satisfy the following conditions:
\begin{enumerate}[(i)]
\item For $p \geq 5$, if ${\underline{v}}: k=p-3, a_{n-p+3} \neq 0, a_j = 0 ~(n-p+4 \leq j \leq n-1), w = 0 = z$ then ${\mathrm{exp}}(G) = p^e, ed(G) \geq 2$ and $||G^{\ast}_0|| = \{ p \}$.  

\item For any odd prime $p$, there exists integers $1 \leq a_{n-p+2}, w, z \leq p-1$ such that the parameters ${\underline{v}}: k=p-2, a_{n-p+2}, a_j = 0 ~(n-p+3 \leq j \leq n-1), w, z$ satisfy ${\mathrm{exp}}(G) = p^e, {\mathrm{ed}}(G) \geq 2$ and $||G^{\ast}_0|| = \{ p^2 \}$. 

\item For any prime $p \geq 5$ (resp. for any odd prime $p$), there exists integers $1 \leq a_{n-p+2}, z \leq p-1$ and $0 \leq w \leq p-1$ such that the parameters ${\underline{v}}: k=p-2, a_{n-p+2}, a_j = 0 ~(n-p+3 \leq j \leq n-1), w, z$ satisfy ${\mathrm{exp}}(G) = p^e, ed(G) \geq 2$ and $||G^{\ast}_0|| = \{ p, p^2 \}$ and $G$ is not an exceptional type (resp. $G$ is of $p$-exceptional type).
\end{enumerate}
\end{theorem}

\noindent {\bf Proof.} (i) We follow the notations of {\ref{mie-struc}}. Here $a_{n-p+3} \neq 0$ ensures $G_1$ is not abelian. Since $G_2$ is abelian of order $p^{(e-1)(p-1)+1}$, using {\ref{abel-sub}} we have ${\mathrm{exp}}(G_2) = p^e$ and ${\mathrm{exp}}(G_3) = p^{e-1}$ and hence $ed(G) \geq 2$. Finally using {\ref{G0lem1}} we have $||G^{\ast}_0|| = \{ p \}$. For the rest of the proof we need the following lemma: 

\begin{lemma}\label{quad-eqn} Let $p$ be an odd prime and consider the equation 
\[
f(X) = a + bX + cX^2
\]
with $a, b, c \in {\mathbb F}_p$ and $c \neq 0$.
\begin{enumerate}
\item There exists $a,b,c \in {\mathbb F}^{\ast}_p$ such that $f(X)=0$ has no solution in ${\mathbb F}_p$. 

\item Given any $b, c \in {\mathbb F}^{\ast}_p$ there exists an $a \in {\mathbb F}_p$ so that $f(X)$ has two distinct solutions.
\end{enumerate}
\end{lemma}

\noindent {\bf Proof.} The equation $f(X)=0$ has a solution if and only if the discriminant $\Delta = b^2 - 4ac \in ({\mathbb F}^{\ast}_p)^2 \cup \{ 0 \}$. While $\Delta = 0$, it has a unique solution and while $\Delta \in ({\mathbb F}^{\ast}_p)^2$ has two distinct solutions.

\noindent (i) For $r^{\prime} \not\in ({\mathbb F}^{\ast}_p)^2 \cup \{ 0 \}$ and $\lambda \in {\mathbb F}^{\ast}_p$ choose $a = 4^{-1} \lambda^2, b = \lambda^2, c = \lambda^2 - r^{\prime} \in {\mathbb F}_p$. Then $\Delta = \lambda^2 r^{\prime} \not\in ({\mathbb F}^{\ast}_p)^2$. (ii) Given $a, b \in {\mathbb F}^{\ast}_p$, choose $r \in ({\mathbb F}^{\ast}_p)^2$ and solve for $r = b^2 - 4ac$ in $c$.  

\subsection*{Proof of {\ref{examples-ed2}} (ii)-(iii)} (ii) First we consider any three such integers $1 \leq a_{n-p+3}, w, z \leq p-1$ so that the conditions are satisfied, and we look forward to the required condition. Everything except $||G^{\ast}_0|| = \{ p^2 \}$ required to be checked. Since the elements $ss^{\zeta}_1 ~(\zeta \in {\mathbb F}_p)$ represents all the $p$ distinct $z$-classes in $G^{\ast}_0$, it is enough to prove each of these elements have order $p^2$. By {\ref{G0lem1}}, this is true if and only if $w + z X + \psi X^2$ has no solution in ${\mathbb F}_p$ where $\psi = a_{n-p+3} \not\equiv 0$ mod $p$. These choices on $w,z,\psi \in {\mathbb F}^{\ast}_p$ can be made by {\ref{quad-eqn}}(i).  

(iii) For the first part, we take $z, \psi = a_{n-p+2} \not\equiv 0$ mod $p$. As in lemma \ref{quad-eqn}(ii), there exists $w \in {\mathbb F}_p$ so that the polynomial $w + z X + \psi X^2$ has exactly two solutions in $\zeta_1, \zeta_2 \in {\mathbb F}_p$, which require $z^2 - 4w\psi \in ({\mathbb F}^{\ast}_p)^2$. By \ref{G0lem1}, the equation $(ss^{\zeta}_1)^p = s^{w + z \zeta + \psi \zeta^2}_{n-1} = 1$ is solved only for the choices $\zeta = \zeta_1, \zeta_2$ in ${\mathbb F}_p$. However, the elements $ss^{\zeta}_1 ~(\zeta \in {\mathbb F}_p)$ represents all the $p$ distinct $z$-classes in $G^{\ast}_0$. Since $p \geq 5$, the remaining $(p-2) \geq 3$ of these representatives has order $p^2$. Proof of the second part is same as above where we need to choose $w = z^2 (4 \psi)^{-1}$ to make sure the polynomial has exactly one solution. This implies $G$ is of $p$-exceptional type. \QED 

\begin{remark}
Note that following the arguments in the proof of \ref{examples-ed2}, for $p = 3$, we obtain a $p$-group of $3^2$-exceptional type. These groups are metabelian (\cite[Thm.5.2]{falc}), and are seen as the only exceptions in the following theorem. Construction of groups of $p^2$-exceptional type for primes $p \geq 5$ require some more serious combinatorics which is not discussed in the paper.    
\end{remark}

\begin{theorem}\label{metab-non-p2} Let $p \geq 5$ be a prime and $G$ be a finite $p$-group of maximal class of order $p^n \geq p^{p+1}$. If $G$ is of $p^2$-exceptional type, then it cannot be metabelian and the nilpotency class of $G_1$ cannot be $\leq 2$.
\end{theorem}

\noindent {\bf Proof.} Let $G$ be metabelian. We choose $w, z, \psi$ mod $p$ as in {\ref{G0lem1}}. If $G$ is of $p^2$-exceptional type, then $w + z \zeta + \psi \zeta^2 \equiv 0$ mod $p$ has $p-1 \geq 4$ solutions in $\zeta$. Since this is a quadratic equation in $\zeta$ we have $w, z, \psi = 0$. But this implies $||G^{\ast}_0|| = \{ p \}$, a contradiction. The second part of the proof follows from \cite[Thm.2.2]{fggs}. \QED

\smallskip
\begin{center}
{ACKNOWLEDGEMENT}
\end{center}

\noindent The author would like to thank J\"{u}rgen M\"{u}ller for the examples of groups of small order using GAP.  

 
\bibliographystyle{plain} 
\bibliography{SpectrumMaxClassV3}

\end{document}